\documentclass[a4paper,11pt,oneside]{article}

\usepackage{amsmath,amsthm,amsfonts,amssymb,amscd}
\usepackage{mathtools}
\usepackage{float}
\usepackage{subfig} 
\usepackage{enumerate}
\usepackage{booktabs}
\usepackage[noend]{algpseudocode}
\usepackage[normalem]{ulem}
\usepackage{graphicx,bm,xcolor}
\usepackage{algorithm}
\usepackage{algpseudocode}
\usepackage[hidelinks]{hyperref}
\usepackage{caption}
\usepackage[T1]{fontenc}
\usepackage{t1enc}
\usepackage[polish,german,english]{babel}
\usepackage{verbatim}
\usepackage{diagbox}
\usepackage{colortbl}
\usepackage{siunitx}

\usepackage{url}
\usepackage{authblk}

\usepackage{pgf,tikz,pgfplots}
\usepackage{mathrsfs}
\usetikzlibrary{arrows}
\usetikzlibrary{shapes.misc}
\usetikzlibrary{positioning}
\usetikzlibrary{decorations.pathreplacing} 

\tikzset{cross/.style={cross out, draw=black, minimum size=2*(#1-\pgflinewidth), inner sep=0pt, outer sep=0pt},
cross/.default={1pt}}

\usepackage{tikz-cd}
\usepackage{svg}

\algnewcommand{\algorithmicgoto}{\textbf{go to}}
\algnewcommand{\Goto}[1]{\algorithmicgoto~\ref{#1}}

\usepackage[colorinlistoftodos,prependcaption,textsize=scriptsize,color=red!40!white]{todonotes}\usepackage{totcount}

\definecolor{orange}{RGB}{230, 159, 0}
\definecolor{skyblue}{RGB}{86, 180, 233}
\definecolor{yellow}{RGB}{240, 228, 66}
\definecolor{blue}{RGB}{0, 114, 178}
\definecolor{vermillion}{RGB}{213, 94, 0}

\definecolor{brewer1}{HTML}{A6CEE3}
\definecolor{brewer2}{HTML}{1F78B4}
\definecolor{brewer3}{HTML}{B2DF8A}
\definecolor{brewer4}{HTML}{33A02C}

\newcommand\auth{\textcolor{black}}

\usepackage[top=2cm, bottom=2cm, left=2cm, right=2cm]{geometry}

\theoremstyle{plain} 
\newtheorem{theorem}{Theorem}[section]

\newtheorem{remark}[theorem]{Remark}

\theoremstyle{definition} %

\theoremstyle{remark} %

\DeclareMathOperator*{\dG}{dG}
\DeclareMathOperator*{\cG}{cG}

\newcommand\restrict[1]{\raisebox{-.5ex}{$|$}_{#1}}

\definecolor{brewer1}{HTML}{A6CEE3}
\definecolor{brewer2}{HTML}{1F78B4}
\definecolor{brewer3}{HTML}{B2DF8A}
\definecolor{brewer4}{HTML}{33A02C}

\begin{document}


\title{
A monolithic space-time temporal multirate finite element framework
for interface and volume coupled problems
}
\author[1,2]{Julian Roth}
\author[3]{Martyna Soszy\'nska}
\author[3]{Thomas Richter}
\author[1,2,4]{Thomas Wick}

\affil[1]{Leibniz Universit\"at Hannover, Institut f\"ur Angewandte
  Mathematik, Welfengarten 1, 30167 Hannover, Germany}
\affil[2]{Universit\'e Paris-Saclay, LMPS - Laboratoire de Mecanique Paris-Saclay,
91190 Gif-sur-Yvette, France}
\affil[3]{Otto-von-Guericke Universit\"at Magdeburg, Universit\"atsplatz 2, 39106 Magdeburg, Germany}
\affil[4]{Cluster of Excellence PhoenixD (Photonics, Optics, and
Engineering - Innovation Across Disciplines), Leibniz University Hannover, Germany}

\date{}
\maketitle

\setcounter{page}{1}

\begin{abstract}
In this work, we propose and computationally investigate a monolithic space-time 
multirate scheme for coupled problems. The novelty lies in the monolithic 
formulation of the multirate approach as this requires a careful design of the functional framework, corresponding discretization, and implementation.
Our method of choice is a tensor-product Galerkin space-time discretization. 
The developments are carried out for both prototype interface- and volume coupled problems such as coupled wave-heat-problems and a displacement equation 
coupled to Darcy flow in a poro-elastic medium. The latter is applied to 
the well-known Mandel's benchmark and a three-dimensional footing problem. Detailed computational investigations 
and convergence analyses give evidence that our monolithic multirate framework 
performs well.
\end{abstract}

\section{Introduction}
\noindent This work is devoted to space-time coupled problems with a
monolithic numerical solution using different temporal meshes, e.g. different time step sizes for each subproblem. Such schemes using different time meshes
for the subproblems are known as multirate approaches.
Multirate schemes are well-known specifically in ordinary differential equations \cite{Logg2003a,Logg2003b,SavHuVe07,SAVCENCO2009323}, porous media 
where Darcy flow and geomechanics couple
\cite{Dean2006,Almani2016,Bause2017,RyMaHeRoh15,ALMANI20192682,Almani2016_phd,RYBAK2014327}, phase-field fracture porous media \cite{JaWheWi21,AlLeeWheWi17,WheWiLee20},
Stokes-Darcy coupling \cite{ShZhLa13}, wave equation and Richardson equation \cite{SOCHALA20092122}, advection equations \cite{SCHLEGEL2009345}, and coupled flow and transport \cite{BrKoeBau22}. In \cite{Soszynska2021,Soszynska_phd}, 
the coupling of heat and wave and thermoelasticity was investigated, and therein
the time meshes are selected using a posteriori error control using
the dual-weighted residual method \cite{BeRa01}.

To the best of our knowledge, in all these studies, partitioned approaches, 
e.g. \cite{Dean2006,JaWheWi21,Almani2016} or one-way coupling techniques \cite{BrKoeBau22} were employed. 
This is reasonable since here multirate iterative 
coupling procedures allow using different programming codes for the subproblems 
and employing multirate schemes can yield very efficient numerical solution 
methods. On the other hand, 
monolithic methods are well-known for their robustness, and if preconditioners 
are available, they can also be more efficient than partitioned approaches. 
Such a monolithic scheme with a focus on the temporal scales 
(not yet multirate though)
for a challenging coupled problem, namely fluid-structure interaction, was proposed in \cite{MaKlWaGe15}. 

The governing model and discretization are based on space-time 
schemes \cite{LaStein19,TeBeLi92,TeBeMiJo92,BaGeRa10,SchVe08,RoThiKoeWi2023,DoeWie22} in which both the temporal and the spatial parts are discretized with 
Galerkin finite elements.
Specifically, in this work we are interested in tensor-product 
space-time finite element discretizations. Here, for efficiency reasons, the problems 
are solved on time-slabs (for the terminology time-slab; see \cite{HuHu88,JOHNSON1993117}), e.g. $(0, T_1), (T_1, T_2), ..., (T_{l-1}, T_l)$, where $T_l:= T$, rather than the entire time interval $(0,T)$.

The key objective of this work is to design and computationally investigate a 
space-time multirate framework in which the coupled system is treated in a 
variational-monolithic fashion. Therein, the coupling can be of interface 
type or achieved via volume coupling. An example for interface coupling 
is fluid-structure interaction \cite{BaKeTe13,Ri17_fsi} and an example 
for volume coupling are porous media equations \cite{Coussy2004} or phase-field fractured porous media \cite{Wick2020PFF}.
Our proposed concepts can deal with both situations. In order to develop and 
analyze our concepts, we however concentrate on prototype settings in the 
interface couplings such as a (linear) heat-wave equation system and 
a linearly coupled poro-elasticity problem, the so-called Biot equations 
\cite{Biot41a,Biot41b,Biot55,Biot7172,To92}.
One of the main tasks is to design functional frameworks and corresponding 
space-time discretizations for such multirate formulations. Here, the key 
is the correct prescription of coupling conditions within the function spaces
and their resulting discrete forms and the corresponding implementation. 
To the best of our knowledge such monolithic space-time multirate frameworks 
are novel. The systems are substantiated with various numerical tests 
for spatial 1d heat-wave equations, spatial 2d heat-wave equations, and the Biot 
system solved on the Mandel benchmark problem and a spatial 3d footing problem. The performances are evaluated 
in terms of goal functionals, their accuracy, and the cost complexity in terms 
of the required space-time degrees of freedom.

The outline of this paper is as follows. In Section \ref{sec_methodology}, 
the monolithic space-time methodology is introduced for abstract coupled problems. 
Next in Section \ref{sec_model_problems}, our two model problems are introduced, namely 
the interface-coupled heat-wave system, and secondly, the volume-coupled poro-elasticity 
problem. In Section \ref{sec_numerical_tests} five numerical tests are conducted.
Finally, in Section \ref{sec_conclusions}, our work is summarized.

\section{Monolithic space-time multirate methodology}
\label{sec_methodology}
Let $W_1,W_2,W_3$ and $V_1,V_2,V_3$ be Banach spaces, and 
$\mathcal{A}_1:W_1\mapsto V_1^*$ and $\mathcal{A}_2:W_2\mapsto V_2^*$ are possibly nonlinear spatio-temporal operators, 
$\mathcal{B}:W_3\mapsto V_3^*$ is the interface coupling operator, where $V_i^*, i=1,2,3$ denote the dual spaces of $V_i,i=1,2,3$. In more detail, the differential operators 
$\mathcal{A}_1,\mathcal{A}_2,\mathcal{B}$ act between time-dependent 
Sobolev spaces.
In the following, we consider the abstract interface coupled multiphysics problem first in strong form, i.e., the operator formulations
\begin{subequations}\label{eq_abstract_interface_problem}
\begin{align}
    \mathcal{A}_1(u_1) &= f_1 \phantom{0}\qquad\quad\text{in } \Omega_1 \times I, \\
    \mathcal{A}_2(u_2) &= f_2 \phantom{0}\qquad\quad\text{in } \Omega_2  \times I, \\
    \mathcal{B}(u_1, u_2) &= 0 \phantom{f_1}\qquad\quad\text{in } \Gamma  \times I,
\end{align}
\end{subequations}
with $u_1: \bar{\Omega}_1 \times \bar{I} \rightarrow \mathbb{R}^{d_1}$, $u_2: \bar{\Omega}_2 \times \bar{I} \rightarrow \mathbb{R}^{d_2}$,
and the abstract volume coupled multiphysics problem
\begin{subequations}\label{eq_abstract_volume_problem}
\begin{align}
    \mathcal{A}_1(u_1, u_2) = f_1  \qquad\quad \text{in } \Omega \times I, \\
    \mathcal{A}_2(u_1, u_2) = f_2 \qquad\quad \text{in }  \Omega \times I,
\end{align}
\end{subequations}
with $u_1: \bar{\Omega} \times \bar{I} \rightarrow \mathbb{R}^{d_1}$, $u_2: \bar{\Omega} \times \bar{I} \rightarrow \mathbb{R}^{d_2}$, and $d_1, d_2 \in \mathbb{N}$ depending on whether the solution fields are scalar-valued ($d_i = 1$) or vector-valued ($d_i > 1$).
Herein, $f_1: {\Omega}_1 \times {I} \rightarrow \mathbb{R}^{d_1}, 
f_2: {\Omega}_2 \times {I} \rightarrow \mathbb{R}^{d_2}$
are sufficiently regular right hand side functions. 
Additionally, these problems need to be completed by suitable initial and boundary conditions. The temporal domain is denoted by $I := (0,T)$ and the spatial domain $\Omega \subset \mathbb{R}^d$ with $d \in \{1,2,3\}$ decomposes for the interface coupled problem into the disjoint subdomains $\Omega_1$ and $\Omega_2$ with common interface $\Gamma := \bar{\Omega}_1 \cap \bar{\Omega}_2$.
\begin{remark}
 The methodology proposed in this paper can also be applied to multiphysics problems that contain both interface and volume coupling.
\end{remark}
Choosing a suitable continuous spatial function space $V(\Omega) = V^1(\Omega_1) \times V^2(\Omega_2)$ for which the interface conditions are well-defined or $V(\Omega) = V^1(\Omega) \times V^2(\Omega)$, we define the space-time function space $X(I, V(\Omega))$ as
\begin{align*}
    X(I, V(\Omega)) := L^2(I, V(\Omega)) \cap H^1(I, V^\ast(\Omega))
\end{align*}
with the dual space of $V(\Omega)$ being denoted as $V^\ast(\Omega) = L(V(\Omega), \mathbb{R})$. 
We briefly notice that with the abstract notation at the beginning, we have 
in this work $W_1\times W_2 = V_1\times V_2 := X(I, V(\Omega))$. 

Thus, for the abstract interface and volume coupled problems, we get a continuous spatio-temporal variational formulation:
    Find $U := \begin{pmatrix}
        u_1 \\ u_2
    \end{pmatrix} \in X(I, V(\Omega))$ such that
    \begin{align}
\label{eq_abstract_form}
        A(U)(\Phi) = F(\Phi) \qquad \forall \Phi := \begin{pmatrix}
        \Phi_1 \\ \Phi_2
        \end{pmatrix} \in \auth{X(I, V(\Omega))}.
    \end{align}
The weak space-time formulations of two model problems for interface and volume coupling will be described in Section \ref{sec_model_problems}.
We notice that $V(\Omega)$ is specified below and we also refer the reader to the literature for more details on the continuous level function spaces for parabolic 
problems \cite{LaStein19} and the wave equation \cite{BaGeRa10}. The latter is 
represented as mixed-in-time formulation in \eqref{eq_abstract_form}
with
\[
X(I, V(\Omega)) := \{U |\; U\in L^2(I,V(\Omega)), \quad \partial_t U \in L^2(I,L^2(\Omega)), \quad \partial_{tt} U \in L^2(I,V(\Omega)^\ast)\}.
\]
The corresponding isometric isomorphic function spaces 
for the displacements and velocities are:
\begin{align*}
    X(I, V^u(\Omega)) & := {L^2(I,V(\Omega)) \cap H^1(I, L^2(\Omega)) \cap H^2\left(I, (V(\Omega))^\ast\right)}, \\
    X(I, V^v(\Omega)) & := {L^2(I, L^2(\Omega)) \cap H^1\left(I, (V(\Omega))^\ast\right)},                                         \\
    X(I, V(\Omega))   & :=  X(I, V^u(\Omega)) \times  X(I, V^v(\Omega)).
\end{align*}
For the porous 
media function spaces all details are given in Section \ref{sec_model_problem_poroelasticity}.

In the following, we will use the notation
\begin{align*}
    (f,g) := (f,g)_{L^2(\Omega)} := \int_\Omega f \cdot g\ \mathrm{d}x, \qquad (\!(f,g)\!) := (f,g)_{L^2(I, L^2(\Omega))} := \int_I (f, g)\ \mathrm{d}t, \\
    \langle f,g \rangle_{{\Gamma}} := \langle f,g \rangle_{L^2(\Gamma)} := \int_\Gamma f \cdot g\ \mathrm{d}s, \qquad (\!\langle f,g\rangle\!)_{{\Gamma \times I}} := (f,g)_{L^2(I, L^2(\Gamma))} := \int_I \langle f, g\rangle_{{\Gamma}} \mathrm{d}t, \\
{    (f,g)_{\Omega_i} := (f,g)_{L^2(\Omega_i)} := \int_{\Omega_i} f \cdot g\ \mathrm{d}x, \qquad (\!(f,g)\!)_{\Omega_i \times I} := (f,g)_{L^2(I, L^2(\Omega_i))} := \int_I (f, g)_{\Omega_i}\ \mathrm{d}t,}
\end{align*}
{for $i \in \lbrace 1,2 \rbrace$.}
In this notation, $f \cdot g$ represents the Euclidean inner product if $f$ and $g$ are scalar- or vector-valued and it represents the Frobenius inner product if $f$ and $g$ are matrices.

\subsection{Discretization in time and space}
\label{sec_discretization_space_time}
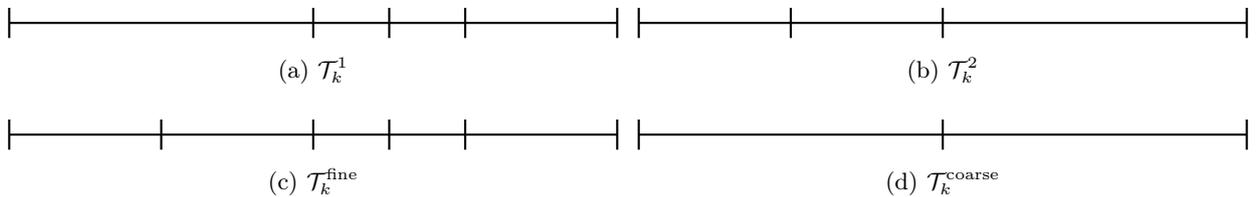
\begin{figure}[H]
     \centering
     \subfloat[$\mathcal{T}_k^1$]{
        \begin{tikzpicture}[scale=1]
            \draw [-, thick] (0,0) -- (8,0);
            \foreach \x in {0,4,5,6,8}
                \draw [thick] (\x,0.2) -- (\x,-0.2);
        \end{tikzpicture}
      }%
    \subfloat[$\mathcal{T}_k^2$]{
        \begin{tikzpicture}[scale=1]
            \draw [-, thick] (0,0) -- (8,0);
            \foreach \x in {0,2,4,8}
                \draw [thick] (\x,0.2) -- (\x,-0.2);
        \end{tikzpicture}
    } \\
    \subfloat[$\mathcal{T}_k^{\text{fine}}$]{
        \begin{tikzpicture}[scale=1]
            \draw [-, thick] (0,0) -- (8,0);
            \foreach \x in {0,2,4,5,6,8}
                \draw [thick] (\x,0.2) -- (\x,-0.2);
        \end{tikzpicture}
      }%
    \subfloat[$\mathcal{T}_k^{\text{coarse}}$]{
        \begin{tikzpicture}[scale=1]
            \draw [-, thick] (0,0) -- (8,0);
            \foreach \x in {0,4,8}
                \draw [thick] (\x,0.2) -- (\x,-0.2);
        \end{tikzpicture}
    }
    \caption{Example for possible temporal meshes $\mathcal{T}_k^1, \mathcal{T}_k^2, \mathcal{T}_k^{\text{fine}}$ and $\mathcal{T}_k^{\text{coarse}}$}
    \label{fig:temporal_mesh_definitions}
\end{figure}

Let the common coarse mesh 
\[
\mathcal{T}_k^{\text{coarse}} := \left\lbrace I_m^{\text{coarse}} := (t_{m-1}^{\text{coarse}}, t_m^{\text{coarse}}) \mid 1 \leq m \leq M^{\text{coarse}} \right\rbrace
\]
be a partitioning of time, i.e.  $\bar{I} = [0,T] = \bigcup_{m = 1}^{M^{\text{coarse}}} \bar{I}_m^{\text{coarse}}$, from which the 
temporal meshes $\mathcal{T}_k^1$ and $\mathcal{T}_k^2$ for the subproblems originate and are being created by local refinement. 
Consequently, 
\[
\mathcal{T}_k^1 := \left\lbrace I_m^1 := (t_{m-1}^1, t_m^1) \mid 1 \leq m \leq M^1 \right\rbrace
\]
and 
\[
\mathcal{T}_k^2 := \left\lbrace I_m^2 := (t_{m-1}^2, t_m^2) \mid 1 \leq m \leq M^2 \right\rbrace
\]
are also partitionings of time for the individual subproblems, i.e.  
\[
\bar{I} = [0,T] = \bigcup_{m = 1}^{M^1} \bar{I}_m^1 = \bigcup_{m = 1}^{M^2} \bar{I}_m^2, 
\]
with $M^i$ being the number of temporal elements for the solution variable $u^i$.
Additionally, we define by 
\begin{align*}
    \mathcal{T}_k^{\text{fine}} := \mathcal{T}_k^1 \cap \mathcal{T}_k^2 := \left\lbrace I_m^1 \cap I_m^2 \mid  I_m^1 \in \mathcal{T}_k^1,\ I_m^2 \in \mathcal{T}_k^2 \right\rbrace    
\end{align*}
the common fine temporal mesh.
These definitions can be visualized by the example in Figure \ref{fig:temporal_mesh_definitions}.
\noindent Then as an intermediate step for discontinuous in time function spaces, we introduce the broken continuous level function space
{
\begin{align*}
    \tilde{X}((\mathcal{T}_k^1, \mathcal{T}_k^2), V(\Omega)) := \left\lbrace \begin{pmatrix}
    u_1 \\ 
    u_2 
    \end{pmatrix} \in L^2(I, L^2(\Omega)^{d_1+d_2})\middle| u_i\restrict{I_{m}^i} \in X(I_{m}^i, V^i(\Omega)),\  I_{m}^i \in \mathcal{T}_k^i,\ i \in \{1,2\} \right\rbrace,
\end{align*}
with the property $X(I, V(\Omega)) \subset \tilde{X}((\mathcal{T}_k^1, \mathcal{T}_k^2), V(\Omega))$.
}
which allows for discontinuities of the solution between temporal elements of $\mathcal{T}_k^1$ and $\mathcal{T}_k^2$.
 By introducing jump terms between temporal elements, we can then define an abstract discontinuous spatio-temporal variational formulation:
    Find $U := \begin{pmatrix}
        u_1 \\ u_2
    \end{pmatrix} \in \tilde{X}((\mathcal{T}_k^1, \mathcal{T}_k^2), V(\Omega)) $ such that
    \begin{align*}
        \tilde{A}(U)(\Phi) = \tilde{F}(\Phi) \qquad \forall \Phi := \begin{pmatrix}
        \Phi_1 \\ \Phi_2
        \end{pmatrix} \in \auth{\tilde{X}((\mathcal{T}_k^1, \mathcal{T}_k^2), V(\Omega))}.
    \end{align*}
We can thus define semi-discrete continuous Galerkin ($\cG$) 
and discontinuous Galerkin ($\dG$) in time function spaces by 
\begin{align*}
    X_k^{\cG(r)}((\mathcal{T}_k^1, \mathcal{T}_k^2), V(\Omega)) := \left\lbrace \begin{pmatrix}
    u_1 \\ 
    u_2 
    \end{pmatrix} \in C(\bar{I}, L^2(\Omega)^{d_1+d_2})\middle| u_i\restrict{I_{m}^i} \in P_r(I_{m}^i, V^i(\Omega)),\  I_{m}^i \in \mathcal{T}_k^i,\ i \in \{1,2\} \right\rbrace ,
\end{align*}
and
\begin{align*}
    X_k^{\dG(r)}((\mathcal{T}_k^1, \mathcal{T}_k^2), V(\Omega)) := \left\lbrace \begin{pmatrix}
    u_1 \\ 
    u_2 
    \end{pmatrix} \in L^2(\bar{I}, L^2(\Omega)^{d_1+d_2})\middle| u_i\restrict{I_{m}^i} \in P_r(I_{m}^i, V^i(\Omega)),\  I_{m}^i \in \mathcal{T}_k^i,\ i \in \{1,2\} \right\rbrace.
\end{align*}
Here, $P_r(I_m, Y)$ is the space of polynomials of order $r$, which map from the time interval $I_m$ into the space $Y$.
Finally, we replace the continuous spatial function space $V(\Omega) = V^1(\Omega) \times V^2(\Omega)$ by conforming finite element subspaces $V_h^1(\mathcal{T}_h^1) \subset V^1(\Omega)$ and $V_h^2(\mathcal{T}_h^2) \subset V^2(\Omega)$, where we define $V_h := V_h(\mathcal{T}_h) := V_h^1(\mathcal{T}_h^1) \times V_h^2(\mathcal{T}_h^2)$. Then the fully discrete function spaces are
\begin{align*}
    X_k^{\cG(r)}\left((\mathcal{T}_k^1, \mathcal{T}_k^2), V_h(\mathcal{T}_h)\right) \subset X_k^{\cG(r)}((\mathcal{T}_k^1, \mathcal{T}_k^2), V(\Omega)),
\end{align*}
and 
\begin{align*}
    X_k^{\dG(r)}\left((\mathcal{T}_k^1, \mathcal{T}_k^2), V_h(\mathcal{T}_h)\right) \subset X_k^{\dG(r)}((\mathcal{T}_k^1, \mathcal{T}_k^2), V(\Omega)).
\end{align*}
Note that in this work we restrict ourselves to spatial meshes $\mathcal{T}_h^1$ and $\mathcal{T}_h^2$ that are constant in time, but with a few minor modifications, this framework can be extended to dynamical spatial meshes which change for each coarse temporal element $I_m \in \mathcal{T}_k^{\text{coarse}}$.

\subsection{Temporal multirate tensor-product space-time FEM}
\label{sec_temporal_multirate_tp_st_fem}
\begin{figure}[H]
    \begin{center}
    \begin{tikzpicture}
    \draw (2,2) -- (2,-2) -- (-2,-2) -- (-2,2) -- (2,2);
    \draw[dashed] (0.5, 2) -- (0.5, -2);
    \draw[dashed] (2, -0.5) -- (-2, -0.5);
    \node at (-0.75,0.75) {$A_1$};
    \node at (-0.75,-1.25) {$B_2$};
    \node at (1.25,-1.25) {$A_2$};
    \node at (1.25,0.75) {$B_1$};
    \draw [decorate,
        decoration = {brace}] (2.3,2) --  (2.3,-2);
      \node[rotate=90] at (2.7,0.0) {test functions};
      \draw [decorate,
        decoration = {brace}] (-2,2.3) --  (2,2.3);
      \node at (0.0,2.7) {trial functions};
    \end{tikzpicture}
    \end{center}
    \caption{We show a snapshot of a block matrix which is solved in every time-step of the coarse mesh $\mathcal{T}_k^{\text{coarse}}$. The off-diagonal entries correspond to the coupling conditions.
    We emphasize that this sketch formally represents both interface-coupling and volume-coupling. In interface-coupling, $B_1$ and $B_2$ represent cross terms associated to the interface trial and test functions, while in volume-coupling, domain cross terms are contained in $B_1$ and $B_2$. Specific examples of $B_1$ and $B_2$ are provided in Section \ref{sec_model_problems}.}
    \label{multiblock}
\end{figure}
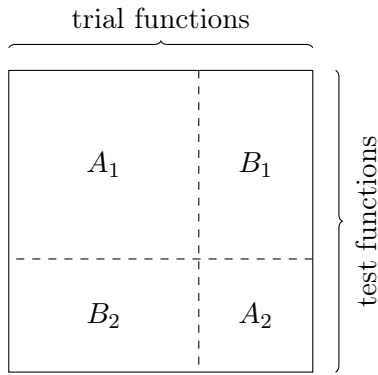

{
Using the different temporal meshes, the space-time discrete variational formulation reads: Find $U_{kh} \in X_k^{\dG(r)}\left((\mathcal{T}_k^1, \mathcal{T}_k^2), V_h(\mathcal{T}_h)\right)$ such that
\begin{align*}
    \tilde{A}(U_{kh})(\Phi_{kh})  = \tilde{F}(\Phi_{kh}) \qquad \forall \Phi_{kh} \in X_k^{\dG(r)}\left((\mathcal{T}_k^1, \mathcal{T}_k^2), V_h(\mathcal{T}_h)\right)
\end{align*}
for a $\dG(r)$ time discretization.
For a continuous-in-time formulation, the discrete variational formulation is: Find $U_{kh} \in X_k^{\cG(r)}\left((\mathcal{T}_k^1, \mathcal{T}_k^2), V_h(\mathcal{T}_h)\right)$ such that
\begin{align*}
    A(U_{kh})(\Phi_{kh})  = F(\Phi_{kh}) \qquad \forall \Phi_{kh} \in X_k^{\dG(r-1)}\left((\mathcal{T}_k^1, \mathcal{T}_k^2), V_h(\mathcal{T}_h)\right),
\end{align*}
but in the following we will only consider discontinuous time discretizations since the FEM discretization for $\cG(r)$ is analogous to the $\dG(r)$ case.
}
Assuming that the PDE is linear, the bilinear form on the left side of this equation can be decomposed into
\begin{equation*}
\tilde{A}(U_{kh})(\Phi_{kh}) = A_1(U_{kh}^1)(\Phi_{kh}^1) + A_2(U_{kh}^2)(\Phi_{kh}^2) + B_1(U_{kh}^2)(\Phi_{kh}^1) + B_2(U_{kh}^1)(\Phi_{kh}^2),
\end{equation*}
where $A_1$ and $A_2$ are bilinear forms for the subproblems and $B_1$ and $B_2$ are bilinear forms containing the coupling terms.
Here, we emphasize that $B_1$ and $B_2$ can include both interface couplings or volume couplings. This should not be confused with the operator $\mathcal{B}$ that has been 
introduced for interface couplings in \eqref{eq_abstract_interface_problem} only, but does appear implicitly in the volume coupled problem \eqref{eq_abstract_volume_problem} as well (cf. Section \ref{sec_model_problems}).
To transition from a standard time marching scheme to this approach, for each time slab one has to build a matrix shown in Figure~\ref{multiblock}. Block $A_1$ represents all the unknowns in the first system without the coupling contributions. $A_2$ is a similar block with all of the unknowns for the second systems. In $B_1$ we have the coupling conditions from the first system and therefore the test functions align with $A_1$ and the trial functions correspond to $A_2$. In $B_2$ we have the coupling conditions of the second system. 
{From a mathematical-numerical viewpoint, it is convenient to work with the space-time form directly by} using tensor-product space-time finite elements as in \cite{BrKoeBau22, RoThiKoeWi2023, FiRoWiChaFau2023}. 
The core idea is that a space-time FEM basis can be created by taking the tensor-product of the spatial and the temporal finite element basis.
\begin{remark}
    Using suitable finite element trial and test spaces, many time marching schemes can be interpreted as variants of Galerkin time discretizations, e.g. the backward Euler method is equivalent to a $\dG(0)$ time 
    discretization under certain assumptions for the quadrature formulae \cite{DeHaTr81,Joh88}.
    Similarly, the Crank-Nicolson method can be interpreted as a variant of the $\cG(1)$ scheme.    
    Moreover, we notice that the Fractional-Step-$\theta$ scheme can be formulated as a Galerkin scheme \cite{MeidnerRichter2014}.
    {Therefore, the tensor-product space-time FEM ansatz includes these time marching schemes as special cases.
    Additionally, we explicitly describe the multirate backward Euler time discretization in \ref{sec_monolithic_multirate_backward_euler} and demonstrate it at an example in \ref{sec_timestepping_mandel}.}
\end{remark}

For efficiency reasons the variational problem does not need to be solved all-at-once on the space-time cylinder $\Omega \times I$, but will be solved forward in time on space-time slabs
\begin{align*}
    S_m := \Omega \times I_m, \qquad I_m \in \mathcal{T}_k^{\text{coarse}}.
\end{align*}
We then have one temporal element for one subproblem and $n \geq 1$ temporal elements for the other subproblem on a given slab $S_m$.
Using this tensor-product space-time FEM discretization, it is fairly straightforward to assemble terms of the variational formulation with trial and test functions from the same subproblem, i.e. the blocks $A_1$ and $A_2$ in Figure \ref{multiblock}.
However, the blocks $B_1$ and $B_2$ with coupling between the subproblems are more complicated, since they require the evaluation of temporal intervals with trial and test functions that belong to different temporal triangulations.
We discuss this problem by considering the interface coupling condition
$\mathcal{B}(u_1, u_2) = u_1 - u_2 = 0$ in $\Gamma  \times I$.
Using penalization to enforce the interface conditions through the weak formulation, we need to be able to assemble terms of the sort $\gamma(\!\langle u_2-u_1,\phi_1\rangle\!)$ \auth{with $\gamma>0$} for the linear system. Note that the only difficulty lies in the assembly of 
$(\!\langle u_2,\phi_1\rangle\!)$.
{
Plugging in the tensor-product space-time finite element ansatz
\begin{align*}
    u_1(x,t) &= \sum_{i_h = 1}^{\#\text{DoF}(\mathcal{T}_h^1)} \sum_{i_k = 1}^{\#\text{DoF}(\mathcal{T}_k^1)} U^1_{i_h, i_k} \phi_h^{1, (i_h)}(x)\phi_k^{1, (i_k)}(t), \\
    u_2(x,t) &= \sum_{i_h = 1}^{\#\text{DoF}(\mathcal{T}_h^2)} \sum_{i_k = 1}^{\#\text{DoF}(\mathcal{T}_k^2)} U^2_{i_h, i_k} \phi_h^{2, (i_h)}(x)\phi_k^{2, (i_k)}(t)
\end{align*}
and assuming that $\phi_1(x,t) = \phi_h^{1, (j_h)}(x) \phi_k^{1,(j_k)}(t)$ for some $j_h, j_k \in \mathbb{N}$, the integral $(\!\langle u_2,\phi_1\rangle\!)$ can be rewritten as
\begin{align}
\label{eq_page_6_formula_1}
    (\!\langle u_2,\phi_1\rangle\!) = \sum_{i_h = 1}^{\#\text{DoF}(\mathcal{T}_h^2)} \sum_{i_k = 1}^{\#\text{DoF}(\mathcal{T}_k^2)} U^2_{i_h, i_k} \cdot \left(\int_\Gamma \phi_h^{2, (i_h)}(x) \phi_h^{1, (j_h)}(x)\ \mathrm{d}x \right) \cdot \left(\int_I \phi_k^{2, (i_k)}(t) \phi_k^{1, (j_k)}(t)\ \mathrm{d}t \right).
\end{align}
}
The spatial integrals over the interface $\Gamma$ can be easily evaluated with most FEM packages. In particular, the evaluation of the interface integral can be found in step-46 of the deal.II \cite{dealii2019design, dealII94} tutorials. Therefore, it only remains to be discussed how the temporal integrals {$\int_I \phi_k^{2, (i_k)}(t) \phi_k^{1, (j_k)}(t)\ \mathrm{d}t$} can be computed, unless the finite element library supports non-matching one dimensional meshes and finite element evaluations.
{
As we have previously mentioned, we assume that the temporal meshes are hierarchical on each slab with
\begin{align*}
    \{I_m \} = \mathcal{T}_k^1 \cap I_m \quad\qquad \text{or} \quad\qquad \{I_m \} = \mathcal{T}_k^2 \cap I_m .
\end{align*}
W.l.o.g. we assume that the first case holds, i.e. $\{I_m \} = \mathcal{T}_k^1 \cap I_m$, and the second temporal mesh is finer than the first mesh.
Then, we can express the temporal basis functions $\phi_k^{1, (j)}$ on the coarser temporal mesh $\mathcal{T}_k^1 \cap I_m$ as linear combinations of the temporal basis functions $\phi_k^{2, (l)}$ on the finer temporal mesh $\mathcal{T}_k^2 \cap I_m$, 
i.e. there exist coefficients $r_{jl} \in \mathbb{R}$ with 
\begin{align*}
    \phi_k^{1, (j)}(t) = \sum_{l = 1}^{\#\text{DoF}(\mathcal{T}_k^2 \cap I_m)} r_{jl} \phi_k^{2, (l)}(t) \qquad \forall t \in I_m.
\end{align*}
Inserting this into the temporal integral in \eqref{eq_page_6_formula_1} from before, we have
\begin{align*}
    \int_{I_m} \phi_k^{2, (i_k)}(t)  \phi_k^{1, (j_k)}(t)\ \mathrm{d}t = \sum_{l = 1}^{\#\text{DoF}(\mathcal{T}_k^2 \cap I_m)} r_{j_k l} \int_{I_m} \phi_k^{2, (i_k)}(t)  \phi_k^{2, (l)}(t)\ \mathrm{d}t.
\end{align*}
This means that the temporal integrals can be computed by evaluating the temporal basis functions on the finer temporal mesh $\mathcal{T}_k^2$ and then restricting it from the left or the right to the coarser temporal mesh $\mathcal{T}_k^1$.
Additionally, we explain this procedure at a concrete example.
}
\auth{We consider} the temporal decompositions $\mathcal{T}_k^1 = \{ (a,b) \}$ and $\mathcal{T}_k^2 = \{ (a,\frac{a+b}{2}), (\frac{a+b}{2},b) \}$ and \auth{$\dG(1)$} basis functions.
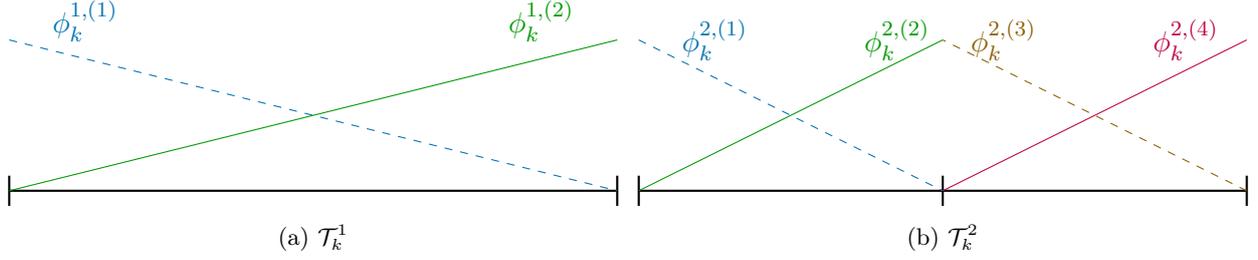
\begin{figure}[H]
     \centering
     \subfloat[$\mathcal{T}_k^1$]{
        \begin{tikzpicture}[scale=1]
            \draw [-, thick] (0,0) -- (8,0);
            \foreach \x in {0,8}
                \draw [thick] (\x,0.2) -- (\x,-0.2);

            \draw [blue, dashed] (0,2) -- (8,0);
            \node [blue] at (1,2.25) {$\phi_k^{1,(1)}$};

            \draw [green!60!black] (0,0) -- (8,2);
            \node [green!60!black] at (7,2.25) {$\phi_k^{1,(2)}$};
        \end{tikzpicture}
      }%
    \subfloat[$\mathcal{T}_k^2$]{
        \begin{tikzpicture}[scale=1]
            \draw [-, thick] (0,0) -- (8,0);
            \foreach \x in {0,4,8}
                \draw [thick] (\x,0.2) -- (\x,-0.2);

            \draw [blue, dashed] (0,2) -- (4,0);
            \node [blue] at (1,2) {$\phi_k^{2,(1)}$};

            \draw [green!60!black] (0,0) -- (4,2);
            \node [green!60!black] at (3.4,2) {$\phi_k^{2,(2)}$};

            \draw [orange!60!black, dashed] (8,0) -- (4,2);
            \node [orange!60!black] at (4.8,2) {$\phi_k^{2,(3)}$};

            \draw [purple] (4,0) -- (8,2);
            \node [purple] at (7.2,2) {$\phi_k^{2,(4)}$};
        \end{tikzpicture}
    }
    \caption{Example of \auth{$\dG(1)$} basis functions on the  temporal meshes $\mathcal{T}_k^1$ and $\mathcal{T}_k^2$}
    \label{fig:dG1_basis_functions}
\end{figure}
\noindent For this example, the matrix corresponding to the term $\auth{(}\!\langle \phi_2,\phi_1\rangle\!\auth{)}$ is thus given by
\begin{align}\label{eq:tensor_product_coupling_matrix}
    \left(M_k R^T\right) \otimes M_h,
\end{align}
where $\otimes$ denotes the Kronecker product, i.e. for $A \in \mathbb{R}^{m \times n}$ and $B \in \mathbb{R}^{p \times q}$ we have
\begin{align*}
    A \otimes B :=  \begin{pmatrix}
        A_{11}B & \cdots & A_{1n}B \\
        \vdots & \ddots & \vdots \\
        A_{m1}B & \cdots & A_{mn}B
    \end{pmatrix} \in \mathbb{R}^{mp \times nq}.
\end{align*}
Furthermore, we have the spatial interface mass matrix
\begin{align*}
    M_h = \left\lbrace \int_\Gamma \phi_h^{2,(j)}(x) \phi_h^{1,(i)}(x)\ \mathrm{d}x 
 \right\rbrace_{1 \leq i \leq \# \text{DoFs}(\mathcal{T}_h^1),\ 1 \leq j \leq \# \text{DoFs}(\mathcal{T}_h^2)},
\end{align*}
the temporal mass matrix on the finer mesh
\begin{align*}
    M_k = \left\lbrace \int_I \phi_k^{2,(j)}(t) \phi_k^{2,(i)}(t)\ \mathrm{d}t 
 \right\rbrace_{1 \leq i,j \leq \auth{4}},
\end{align*}
and 
\begin{align*}
    R = \begin{pmatrix}
        1 & \frac{1}{2} & \auth{\frac{1}{2}} & 0 \\
        0 & \frac{1}{2} & \auth{\frac{1}{2}} & 1
    \end{pmatrix},
\end{align*}
which is the restriction matrix from the fine temporal basis $\Phi_k^{2}$ to the coarse temporal basis $\Phi_k^{1}$, i.e.
\begin{align*}
    \Phi_k^1 = R\Phi_k^2.
\end{align*}
This procedure can then also be extended to other temporal discretizations $\mathcal{T}_k^1$ and $\mathcal{T}_k^2$, as well as other types of $\cG(r)$ or $\dG(r)$ time discretizations, and is possibly already contained in most FEM libraries with multigrid capabilities. 
\begin{remark}\label{rem:tensor_product_coupling_matrix}
    Although the interface coupling matrix in (\ref{eq:tensor_product_coupling_matrix}) has been depicted as a tensor-product of a spatial and a temporal matrix, this is not actually performed in our implementation. Instead, the local space-time matrices are being assembled with the finer temporal basis and the restriction matrix is being applied prior to the distribution of the local contributions to the global matrix. This could possibly be applicable to a matrix-free space-time implementation as well and is also required for e.g. nonlinear problems when the system matrix cannot be formulated as the linear combination of tensor-products of temporal and spatial finite element matrices.
\end{remark}
\begin{remark}
    \auth{
    Since we are using finite elements in time, we can evaluate the solution at any time $t \in I_m$ by evaluating the temporal basis functions at this given time $t$. In particular, this allows the extension of this procedure to nonlinear problems, where after linearization with a Newton scheme one needs to assemble linear problems with temporally varying coefficients that depend on the current iterate.
    }
\end{remark}
\begin{remark}
    We also notice that a conceptionally similar approach has been proposed for ODEs 
    in \cite{Logg2003a,Logg2003b,SavHuVe07}
    and some single non-coupled PDEs \cite{JaLo08}. Therein, the temporal 
    basis is not obtained by restriction
    (as we do), but using the integration
    quadrature points.
\end{remark}

The method provides flexibility and allows for discretization of different terms using separate time methods. Although the size of the system significantly increases, one does not require any subcycling to resolve coupling conditions. Moreover, the resulting scheme is fully implicit which is highly desirable for stiff systems of equations \cite{aiken1985stiff, Hairer1996}. We solve the space-time linear system monolithically with the direct solver UMFPACK \cite{UMFPACK}, but the monolithic system could also be decoupled on the solver level and solved as a partitioned scheme with some iterative solver or some multigrid method; e.g., \cite{Saad2003}.

\section{Model problems in strong forms and weak space-time forms}
\label{sec_model_problems}
To prepare the numerical tests considered in Section \ref{sec_numerical_tests}, we formulate two abstract model problems and derive their weak space-time formulations. 
First, we work with an interface coupled heat and wave equation as in \cite{Soszynska2021}.
Next, we derive the space-time weak form for poroelasticity similar to \cite{Girault2011, Bause2017}.

\subsection{Interface coupling: heat and wave equation}
\label{sec_model_problem_heatwave}
We couple the heat and wave equation across a common interface, which serves as a prototypical example for fluid-structure interaction since we couple a parabolic and a hyperbolic problem across an interface.
Consequently, we name the different subdomains `fluid' and `solid'.
\subsubsection{Strong form}
The wave equation is being solved in the solid domain $\Omega_s$ by finding displacement $u_s: \bar{\Omega}_s \times \bar{I} \rightarrow \mathbb{R}$ and velocity $v_s: \bar{\Omega}_s \times \bar{I} \rightarrow \mathbb{R}$ such that
\begin{align*}
    \partial_t v_s - \lambda \Delta_x u_s - \delta \Delta_x v_s  &= g_s \qquad\quad \text{in } \Omega_s \times I, \\
    \partial_t u_s &= v_s \qquad\quad \text{in } \Omega_s \times I, \\
    u_s = v_s &= 0 \qquad\quad \text{on } \Omega_s \times \{0\}, \\
    u_s = v_s &= 0 \qquad\quad \text{on } \Gamma_D^s \times I, \\
    \lambda \partial_{n_s} u_s + \delta \partial_{n_s} v_s &= 0 \qquad\quad \text{on } \Gamma_N^s \times I,
\end{align*}
where $\lambda, \delta \in \mathbb{R}^+$
and $g_s: \bar{\Omega}_s \times I \rightarrow \mathbb{R}$ is a sufficiently smooth right hand side.
For the (strong) damping term $\delta \Delta_x v_s$ and its mathematical influence, we refer the reader to \cite{GaSqu06}.
The heat equation is being solved in the fluid domain $\Omega_f$ by finding velocity $v_f: \bar{\Omega}_f \times \bar{I} \rightarrow \mathbb{R}$ such that
\begin{align*}
    \partial_t v_f - \nu \Delta_x v_f + \beta \cdot \nabla_x v_f &= g_f \qquad\quad \text{in } \Omega_f \times I,\\
    v_f &= 0 \qquad\quad \text{on } \Omega_f \times \{0\}, \\
    v_f &= 0 \qquad\quad \text{on } \Gamma_D^f \times I, \\
    \partial_{n_f} v_f &= 0 \qquad\quad \text{on } \Gamma_N^f \times I,
\end{align*}
with $\nu \in \mathbb{R}^+$, $\beta \in \mathbb{R}^d$ and $g_f: \bar{\Omega}_f \times I \rightarrow \mathbb{R}$ is a sufficiently smooth right hand side.
Similar to the mesh motion PDE in the Arbitrary Lagrangian Eulerian formulation of fluid-structure interaction, we harmonically extend the solid deformation $u_s$ to the fluid domain by finding the displacement $u_f: \bar{\Omega}_f \times \bar{I} \rightarrow \mathbb{R}$ such that
\begin{align*}
    -\Delta_x u_f &= 0 \qquad\quad \text{in } \Omega_f \times I, \\
    u_f &= 0 \qquad\quad \text{on } \Gamma_D^f \times I, \\
    \partial_{n_f} u_f &= 0 \qquad\quad \text{on } \Gamma_N^f \times I.
\end{align*}
Furthermore, the interface conditions are given by
\begin{align*}
    \lambda \partial_{n_s}u_s + \delta \partial_{n_s}v_s + \nu \partial_{n_f}v_f &= 0 \qquad\quad \text{on } \Gamma \times I,\\
    u_f&= u_s \qquad\quad \text{on } \Gamma \times I, \\
    v_f&= v_s \qquad\quad \text{on } \Gamma \times I.
\end{align*}
In the differential equations above, we assume homogeneous Dirichlet boundary conditions on $\Gamma_D^f$ resp. $\Gamma_D^s$ and homogeneous Neumann boundary conditions on $\Gamma_N^f := \partial \Omega \setminus (\Gamma_D^f \cup \Gamma)$ resp. $\Gamma_N^s := \partial \Omega \setminus (\Gamma_D^s \cup \Gamma)$.
Looking back at the notation in (\ref{eq_abstract_interface_problem}), we have $u_1 := U_f := (u_f, v_f)$, $d_1 = 2$ and otherwise replace the index $1$ by the letter 'f' for the fluid domain, i.e. $\Box_1 := \Box_f$.
Similarly, we have $u_2 := U_s := (u_s, v_s)$, $d_2 = 2$ and otherwise replace the index $2$ by the letter 's' for the solid domain, i.e. $\Box_2 := \Box_s$.
For spatio-temporal operators we have
\begin{align*}
    \mathcal{A}_1(u_1) &:= 
    \begin{pmatrix}
        \partial_t v_f - \nu \Delta_x v_f + \beta \cdot \nabla_x v_f \\
        -\Delta_x u_f
    \end{pmatrix}, \\
    \mathcal{A}_2(u_2) &:= 
    \begin{pmatrix}
        \partial_t v_s - \lambda \Delta_x u_s - \delta \Delta_x v_s \\
        \partial_t u_s - v_s
    \end{pmatrix}, 
\end{align*}
for the right hand sides we have
\begin{align*}
    f_1 &:= \begin{pmatrix}
        g_f \\ 0
    \end{pmatrix}, \\
    f_2 &:= \begin{pmatrix}
        g_s \\ 0
    \end{pmatrix},
\end{align*}
and for the Dirichlet and Neumann interface coupling operators we have
\begin{align*}
    \mathcal{B}_{D,u}(u_1, u_2) &:= u_f - u_s, \\
    \mathcal{B}_{D,v}(u_1, u_2) &:= v_f - v_s, \\
    \mathcal{B}_{N}(u_1, u_2) &:= \lambda \partial_{n_s}u_s + \delta \partial_{n_s}v_s + \nu \partial_{n_f}v_f.
\end{align*}
Herein, $n_f$ is the fluid normal vector, $n_s$ is the solid normal vector and $\partial_n g:= \nabla_x g \cdot n$ is the normal derivative of a function $g$.

The spatial function spaces \cite{Soszynska2021} are given by
\begin{align*}
   V_f(\Omega_f) &= \left(H^1_{0, \Gamma_D^f}(\Omega_f)\right)^2, \\
   V_s(\Omega_s) &= \left(H^1_{0, \Gamma_D^s}(\Omega_s)\right)^2, \\
   V(\Omega) &= V_f(\Omega_f) \times V_s(\Omega_s),
\end{align*}
i.e. they consist of one time weakly differentiable functions that vanish on the Dirichlet boundary $\Gamma_D^f$ resp. $\Gamma_D^s$.

\subsubsection{Weak space-time form}
Integrating the strong formulation over time and multiplying by a test function, we get the discontinuous in time weak formulation:

Find $U \in \tilde{X}\left((\mathcal{T}_k^f, \mathcal{T}_k^s), V(\Omega)\right)$ such that
\begin{align}\label{eq:variational_form_heat_wave}
    \tilde{A}(U)(\Phi)  = \tilde{F}(\Phi) \qquad \forall \Phi \in \auth{\tilde{X}\left((\mathcal{T}_k^f, \mathcal{T}_k^s), V(\Omega)\right)}.
\end{align}
The bilinear form and right hand side read
\begin{align*}
    \tilde{A}(U)(\Phi) &= {A_1(U_f)(\Phi^f) + A_2(U_s)(\Phi^s) + B_1(U_s)(\Phi^f) + B_2(U_f)(\Phi^s),} \\ 
    A_1(U_f)(\Phi^f) &= \sum_{I_m \in \mathcal{T}_k^{\text{fine}}} \int_{I_m} (\partial_t v_f,\phi^{v_f})_{\Omega_f} + \nu ( \nabla_x v_f,\nabla_x \phi^{v_f})_{\Omega_f} + ( \beta \cdot \nabla_x v_f, \phi^{v_f})_{\Omega_f} + ( \nabla_x u_f,\nabla_x \phi^{u_f})_{\Omega_f} \ \mathrm{d}t  \\
     &\quad + \sum_{I_m \in \mathcal{T}_k^{\text{fine}}} \int_{I_m} -\nu\langle\partial_{n_f} v_f,\phi^{v_f}\rangle_{\Gamma} -\langle\partial_{n_f} u_f,\phi^{u_f}\rangle_{\Gamma} + \frac{\gamma \nu}{h} \langle v_f, \phi^{v_f}\rangle_{\Gamma} + \frac{\gamma}{h} \langle u_f, \phi^{u_f}\rangle_{\Gamma} \ \mathrm{d}t \\
     &\quad + \sum_{m_f = 1}^{M^f-1}([v_f]_{m_f}, \phi_{m_f}^{v_f,+})_{\Omega_f} + (v_{f,0}^{+},\phi_0^{v_f,+})_{\Omega_f},\\
    A_2(U_s)(\Phi^s) &= \sum_{I_m \in \mathcal{T}_k^{\text{fine}}} \int_{I_m} (\partial_t v_s,\phi^{v_s})_{\Omega_s}  + \lambda ( \nabla_x u_s,\nabla_x \phi^{v_s})_{\Omega_s} + \delta ( \nabla_x v_s,\nabla_x \phi^{v_s})_{\Omega_s} + (\partial_t u_s,\phi^{u_s})_{\Omega_s} - (v_s, \phi^{u_s})_{\Omega_s}  \ \mathrm{d}t \\
    &\quad - \sum_{I_m \in \mathcal{T}_k^{\text{fine}}} \int_{I_m} \delta \langle \partial_{n_s} v_s, \phi^{v_s}\rangle_{\Gamma}  \ \mathrm{d}t \\
    &\quad + \sum_{m_s = 1}^{M^s-1}([v_s]_{m_s}, \phi_{m_s}^{v_s,+})_{\Omega_s} + ([u_s]_{m_s}, \phi_{m_s}^{u_s,+})_{\Omega_s} + (v_{s,0}^{+},\phi_0^{v_s,+})_{\Omega_s} + (u_{s,0}^{+},\phi_0^{u_s,+})_{\Omega_s}, \\
    B_1(U_s)(\Phi^f) &= \sum_{I_m \in \mathcal{T}_k^{\text{fine}}} \int_{I_m} -\frac{\gamma \nu}{h} \langle v_s, \phi^{v_f}\rangle_{\Gamma} - \frac{\gamma}{h} \langle u_s, \phi^{u_f}\rangle_{\Gamma} \ \mathrm{d}t,\\
    B_2(U_f)(\Phi^s) &= \sum_{I_m \in \mathcal{T}_k^{\text{fine}}} \int_{I_m} \nu \langle \partial_{n_f} v_f, \phi^{v_s}\rangle_{\Gamma}  \ \mathrm{d}t
\end{align*}
and 
\begin{align*}
    \tilde{F}(\Phi) &= (\!(g_f,\phi^{v_f})\!)_{\Omega_f \times I} + (v_f^0, \phi_0^{v_f,+})_{\Omega_f} + (\!(g_s,\phi^{v_s})\!)_{\Omega_s \times I} + (v_s^0, \phi_0^{v_s,+})_{\Omega_s} + (u_s^0, \phi_0^{u_s,+})_{\Omega_s},
\end{align*}
where
\begin{align*}
    U := (u_f, v_f, u_s, v_s), \quad \Phi := (\phi^{u_f}, \phi^{v_f}, \phi^{u_s}, \phi^{v_s}).
\end{align*}

\subsection{Volume coupling: poroelasticity}
\label{sec_model_problem_poroelasticity}

\subsubsection{Strong form}
The governing equations for poroelasticity read \cite{Coussy2004,LeSchref99}: Find pressure $p: \bar{\Omega} \times \bar{I} \rightarrow \mathbb{R}$ and displacement $u: \bar{\Omega} \times \bar{I} \rightarrow \mathbb{R}^d$ such that
\begin{align*}
    \partial_t(cp + \alpha(\nabla_x \cdot u)) - \frac{1}{\nu}\nabla_x \cdot (K(\nabla_x p - \rho g)) &= q \qquad\quad \text{in } \Omega \times I, \\
    -\nabla_x \cdot \sigma(u) + \alpha \nabla_x p &= f \qquad\quad \text{in } \Omega \times I,
\end{align*}
with the stress tensor
\begin{align*}
    \sigma(u) := \mu(\nabla_x u + (\nabla_x u)^T) + \lambda (\nabla_x \cdot u)I.
\end{align*}
A rigorous mathematical analysis of this problem from poroelasticity can be found in \cite{Showalter2000},
and this coupled system of equations is also known as Biot system \cite{Biot41a,Biot41b,Biot55,Biot7172,To92}.

Framing this in the notation in (\ref{eq_abstract_volume_problem}), we have $u_1 := u$, $d_1 := d$ and otherwise replace the index $1$ by the letter 'u' for the displacement u, i.e. $\Box_1 := \Box_u$.
In the same fashion, we have $u_2 := p$, $d_2 = 1$ and otherwise replace the index 2 by the letter 'p' for the pressure, i.e. $\Box_2 := \Box_p$. For the spatio-temporal operators we have 
\begin{align*}
    \mathcal{A}_1(u_1, u_2) &:= \partial_t(cp + \alpha(\nabla_x \cdot u)) - \frac{1}{\nu}\nabla_x \cdot (K(\nabla_x p - \rho g)),\\
    \mathcal{A}_2(u_1, u_2) &:= -\nabla_x \cdot \sigma(u) + \alpha \nabla_x p,
\end{align*}
and for the right hand sides we have
\begin{align*}
    f_1 &:= q, \\
    f_2 &:= f.
\end{align*}
Assuming homogeneous Dirichlet boundary conditions for the displacement on $\Gamma_D$ and inhomogeneous Neumann/traction boundary conditions 
\begin{align*}
    \sigma(u) \cdot n  &= t
\end{align*}
on $\Gamma_N = \partial\Omega\setminus \Gamma_D$, the spatial function spaces are given by
\begin{align*}
   V_u(\Omega) &= \left(H^1_{0, \Gamma_D}(\Omega)\right)^d, \\
   V_p(\Omega) &= H^1(\Omega), \\
   V(\Omega) &= V_u(\Omega) \times V_p(\Omega).
\end{align*}
{We notice that according to the boundary conditions adopted in Section 
\ref{sec_mandel_problem} and Section \ref{sec_footing_problem}, these function 
spaces need to be refined by implementing further Dirichlet conditions on certain
parts, while Neumann conditions arise naturally in the weak forms as usual.}

\subsubsection{Weak space-time form}
We can now derive the space-time variational formulation for this problem. By integration by parts we get the variational formulation: Find $U := \{u, p\} \in X\left(I, V(\Omega)\right)$ such that
\begin{align*}
    c(\!(\partial_t p, \phi^p)\!) + \alpha (\!(\partial_t (\nabla_x \cdot u), \phi^p)\!) + \frac{K}{\nu}(\!(\nabla_x p, \nabla_x\phi^p)\!) &= (\!(q, \phi^p)\!) + \frac{K\rho}{\nu}(\!(g, \nabla_x\phi^p)\!), \\
    (\!(\mu(\nabla_x u + (\nabla_x u)^T) + \lambda (\nabla_x \cdot u)I, \nabla_x\phi^u)\!) - \alpha(\!(pI, \nabla_x\phi^u)\!) + \alpha (\!\langle pn, \phi^u\rangle\!)_{ \Gamma_{\text{top}} \times I} &= (\!(f, \phi^u)\!)  + (\!\langle t, \phi^u \rangle\!)_{ \Gamma_{\text{top}} \times I}\\
    \forall \Phi := \begin{pmatrix}
    \phi^u \\ 
    \phi^p 
    \end{pmatrix}  \in \auth{X\left(I, V(\Omega)\right)}.
\end{align*} 
Accounting for the discontinuities in time, we thus need to solve the problem: 

Find $U \in \tilde{X}\left((\mathcal{T}_k^u, \mathcal{T}_k^p), V(\Omega)\right)$ such that
\begin{align}\label{eq:variational_form_mandel}
    \tilde{A}(U)(\Phi)  = \tilde{F}(\Phi) \qquad \forall \Phi \in \auth{\tilde{X}\left((\mathcal{T}_k^u, \mathcal{T}_k^p), V(\Omega)\right)}.
\end{align}
The bilinear form and right hand side read
\begin{align*}
    \tilde{A}(U)(\Phi) &= {A_1(u)(\phi^u) + A_2(p)(\phi^p) + B_1(p)(\phi^u) + B_2(u)(\phi^p),} \\ 
    A_1(u)(\phi^u) &= \sum_{I_m \in \mathcal{T}_k^{\text{fine}}} \int_{I_m} (\mu(\nabla_x u + (\nabla_x u)^T) + \lambda (\nabla_x \cdot u)I, \nabla_x\phi^u)\ \mathrm{d}t, \\
    A_2(p)(\phi^p) &= \sum_{I_m \in \mathcal{T}_k^{\text{fine}}} \int_{I_m}  c(\partial_t p, \phi^p) + \frac{K}{\nu}(\nabla_x p, \nabla_x\phi^p) \ \mathrm{d}t + \sum_{m_p = 1}^{M^p-1}c([p]_{m_p}, \phi_{m_p}^{p,+}) + c (p_{0}^{+},\phi_0^{p,+}), \\
    B_1(p)(\phi^u) &= \sum_{I_m \in \mathcal{T}_k^{\text{fine}}} \int_{I_m} - \alpha(pI, \nabla_x\phi^u) + \alpha\langle pn, \phi^u\rangle_{ \Gamma_{\text{top}}} \ \mathrm{d}t, \\
    B_2(u)(\phi^p) &= \sum_{I_m \in \mathcal{T}_k^{\text{fine}}} \int_{I_m} \alpha (\partial_t (\nabla_x \cdot u), \phi^p)\ \mathrm{d}t + \sum_{m_u = 1}^{M^u-1}\alpha([\nabla_x \cdot u]_{m_u}, \phi_{m_u}^{p,+}) + \alpha (\nabla_x \cdot u_{0}^{+},\phi_0^{p,+}).
\end{align*}
and 
\begin{align*}
    \tilde{F}(\Phi) &:= (\!(q, \phi^p)\!) + \frac{K\rho}{\nu}(\!(g, \nabla_x\phi^p)\!) + (\!(f, \phi^u)\!)  + (\!\langle t, \phi^u\rangle\!)_{ \Gamma_{\text{top}} \times I} + \alpha(\nabla_x \cdot u^0, \phi_0^{p,+}) + c(p^0, \phi_0^{p,+}),
\end{align*}
where
\begin{align*}
    U := (u, p), \quad \Phi := (\phi^u, \phi^p).
\end{align*}
A derivation of a multirate time stepping scheme for this space-time formulation using Section \ref{sec_temporal_multirate_tp_st_fem} can be found in \ref{sec_timestepping_mandel}. 

\section{Numerical tests}
\label{sec_numerical_tests}
For verification of our multirate framework, we perform numerical tests on {five} 
different coupled problems.
For the first {three} numerical tests, we perform computations for interface coupled heat and wave equations in 1+1D and 2+1D. For the 1+1D problem, we construct a manufactured solution and for the 2+1D problem, we use Configuration{s} 2.1 {and 2.2}  from \cite{Soszynska2021} with a transport term in the heat equation and a damping term in the wave equation.
As the \auth{fourth} numerical test, we consider Mandel's problem 
\cite{Mandel1953,Cr63,Cheng88,AbChCuDeRo96,DuiMiWi22}, 
a poroelasticity benchmark problem in 2+1D, as an example of volume coupled problems. {In the fifth numerical test, we extend to a 3+1D poroelasticity footing problem \cite{Gaspar2008}.}

The time marching codes (Section \ref{sec_1d_heat_wave} only) have been implemented in FEniCS \cite{fenics2015} and the space-time FEM codes have been written in deal.II \cite{dealii2019design, dealII94}.

\subsection{1+1D heat and wave equation}
\label{sec_1d_heat_wave}
\begin{figure}[H]
     \centering
     \subfloat[Displacement]{
    \includegraphics[width=8cm]{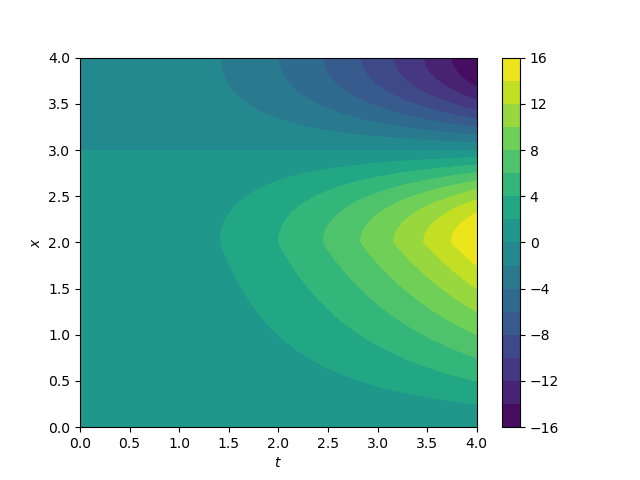}}%
     \subfloat[Velocity]{
    \includegraphics[width=8cm]{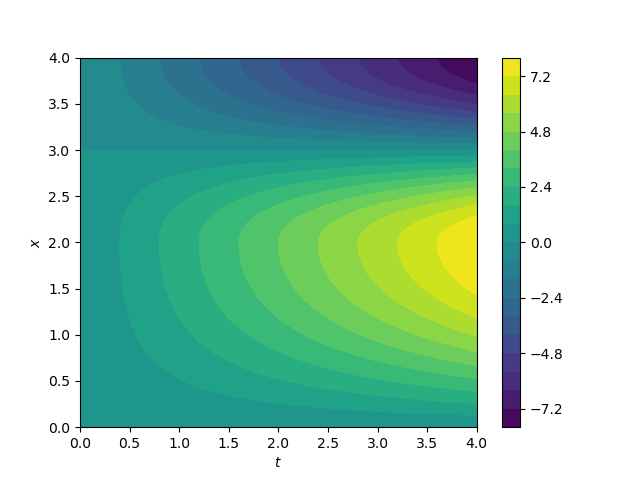}}%
    \caption{Analytical solution for the 1+1D heat and wave equation problem (solid in upper half and fluid in lower half of the images)}
    \label{fig:manufactured_solution_1d}
\end{figure}

As the first numerical test, we consider a one-dimensional problem where we couple the heat and wave equation, cf. Section \ref{sec_model_problem_heatwave}. The time interval is $I:=(0,4)$, the spatial domain for the heat equation - the fluid domain - is $\Omega_f := (0,2)$ and the spatial domain for the wave equation - the solid domain - is $\Omega_s := (2,4)$. The interface is given by $\Gamma := \bar{\Omega}_f \cap \bar{\Omega}_s = \{2\}.$ 
To further simplify the model problem from Section \ref{sec_model_problem_heatwave}, we leave out the transport term in the heat equation and the damping term in the wave equation, i.e. $\beta = 0$ and $\delta = 0$. 
For our test we choose the parameters in the PDEs to be $\nu = 0.001$ and $\lambda = 1000$. The penalty parameter $\gamma$ for the Dirichlet interface conditions is $1000$.
The boundaries are defined as $\Gamma_D^f = \{0\}$, $\Gamma_N^s = \{4\}$ and the variational formulation is shown in (\ref{eq:variational_form_heat_wave}).
We prescribe the analytical solution as
\begin{align*}
    u_s(x,t) &= t^2 \cdot  \cos \left( \frac{\pi(x-2)}{2} \right)  &&\forall (x,t) \in \Omega_s \times I, \\
    u_f(x,t) &= t^2 \cdot \frac{x}{2}  &&\forall (x,t) \in \Omega_f \times I,
\end{align*}
and the velocities
\begin{align*}
    v_s(x,t) &= 2t \cdot  \cos \left( \frac{\pi(x-2)}{2} \right)  &&\forall (x,t) \in \Omega_s \times I, \\
    v_f(x,t) &= 2t \cdot \sin \left( \frac{\pi x}{4} \right)  &&\forall (x,t) \in \Omega_f \times I,
\end{align*}
which are shown in Figure \ref{fig:manufactured_solution_1d}.
\noindent The right hand side functions are
\begin{align*}
    g_s(x,t) &=  2 \cos \left( \frac{\pi(x-2)}{2} \right) +  \frac{\pi^2t^2\lambda\cos(\frac{\pi(x-2)}{2})}{4} &&\forall (x,t) \in  \Omega_s \times I, \\
    g_f(x,t) &= 2 \sin \left( \frac{\pi x}{4} \right) + \frac{\pi^2t\nu\sin(\frac{\pi x}{4})}{8}  &&\forall (x,t) \in \Omega_f \times I.
\end{align*}

\subsubsection{Monolithic multirate \auth{backward Euler} results}
To compare the discrete solutions with the analytic ones, we define 
\begin{align*}
    J(U - U_{kh}) &= \|U - U_{kh}\|_{L^2(I, L^2(\Omega))} { = \sqrt{\|U_f - U_{f,kh}\|_{L^2(I, L^2(\Omega_f))}^2 + \|U_s - U_{s, kh}\|_{L^2(I, L^2(\Omega_s))}^2}}.
\end{align*}
We similarly define specific contributions as 
\begin{align*} 
    J_f(U - U_{kh}) &= \|U_f - U_{f,kh}\|_{L^2(I, L^2(\Omega_f))}, \\
    J_s(U - U_{kh}) &= \|U_s - U_{s, kh}\|_{L^2(I, L^2(\Omega_s))}.
\end{align*}
In the following tables, we will analyze the errors denoted by 
\begin{equation*}
\eta_{kh} \coloneqq J(U - U_{kh}), \hspace*{0.5 cm}\eta_{kh}^f \coloneqq J_f(U - U_{kh}), \hspace*{0.5 cm}\eta_{kh}^s \coloneqq J_s(U - U_{kh}).
\end{equation*}

\begin{table}[H]
  \begin{center}
    \begin{tabular}{cccc|ccc|c}
        \toprule
        $|\mathcal{T}_k^{\text{coarse}}|$ & $|\mathcal{T}_k^f|$ & $|\mathcal{T}_k^s|$ & $|\mathcal{T}_k^f| : |\mathcal{T}_k^s|$ & $\eta_{kh}^f$ & $\eta^s_{kh}$ & $\eta_{kh}$ & {EOC} \\       
        \midrule
      25 & 25 & 25 & 1:1& $1.78 \cdot 10^{-2}$ & $3.73 \cdot 10^{-2}$ & $4.14 \cdot 10^{-2}$ & {-} \\
      50 & 50 & 50 & 1:1 & $8.78 \cdot 10^{-3}$  & $1.90 \cdot 10^{-2}$ &  $2.10 \cdot 10^{-2}$ & {0.98} \\
      100 & 100 & 100 & 1:1 & $4.46 \cdot 10^{-3}$& $9.77 \cdot 10^{-3}$ & $1.08 \cdot 10^{-2}$ & {0.96} \\
      200 & 200 & 200 & 1:1 & $2.57 \cdot 10^{-3}$& $5.16 \cdot 10^{-3}$ & $5.76 \cdot 10^{-3}$ & {0.91} \\
      400 & 400 & 400 & 1:1 & $1.93 \cdot 10^{-3}$  & $2.96 \cdot 10^{-3}$ &   $3.52 \cdot 10^{-3}$ & {0.71} \\
      \bottomrule
    \end{tabular}
    \caption{Errors for fully uniform time-stepping.}
    \label{uniform_1D}
  \end{center}
\end{table}

\begin{table}[H]
  \begin{center}
    \begin{tabular}{cccc|ccc|c}
        \toprule
        $|\mathcal{T}_k^{\text{coarse}}|$ & $|\mathcal{T}_k^f|$ & $|\mathcal{T}_k^s|$ & $|\mathcal{T}_k^f| : |\mathcal{T}_k^s|$ & $\eta_{kh}^f$ & $\eta^s_{kh}$ & $\eta_{kh}$ & {EOC} \\       
        \midrule
      25 & 50 & 25 & 2:1& $2.40 \cdot 10^{-2}$ & $3.86 \cdot 10^{-2}$ & $4.54 \cdot 10^{-2}$ & {-} \\
      50 & 100 & 50 & 2:1 & $1.14 \cdot 10^{-2}$  & $1.98 \cdot 10^{-2}$ &  $2.29 \cdot 10^{-2}$ & {0.99} \\
      100 & 200 & 100 & 2:1 & $5.45 \cdot 10^{-3}$& $1.03 \cdot 10^{-2}$ & $1.17 \cdot 10^{-2}$ & {0.97} \\
      200 & 400 & 200 & 2:1 & $2.72 \cdot 10^{-3}$& $5.53 \cdot 10^{-3}$ & $6.16 \cdot 10^{-3}$ & {0.93} \\
      \bottomrule
    \end{tabular}
    \caption{Errors for one refinement in the fluid domain.}
    \label{refinement_fluid_1D}
  \end{center}
\end{table}

\begin{table}[H]
  \begin{center}
    \begin{tabular}{cccc|ccc|c}
        \toprule
        $|\mathcal{T}_k^{\text{coarse}}|$ & $|\mathcal{T}_k^f|$ & $|\mathcal{T}_k^s|$ & $|\mathcal{T}_k^f| : |\mathcal{T}_k^s|$ & $\eta_{kh}^f$ & $\eta^s_{kh}$ & $\eta_{kh}$ & {EOC} \\       
        \midrule
      25 & 25 & 50 & 1:2& $8.98 \cdot 10^{-3}$ & $2.10 \cdot 10^{-2}$ & $2.28 \cdot 10^{-2}$ & {-} \\
      50 & 50 & 100 & 1:2 & $5.06 \cdot 10^{-3}$  & $1.04 \cdot 10^{-2}$ &  $1.16 \cdot 10^{-2}$ & {0.97} \\
      100 & 100 & 200 & 1:2 & $3.33 \cdot 10^{-3}$& $5.21 \cdot 10^{-3}$ & $6.18 \cdot 10^{-3}$ & {0.91} \\
      200 & 200 & 400 & 1:2 & $2.65 \cdot 10^{-3}$& $2.80 \cdot 10^{-3}$ & $3.86 \cdot 10^{-3}$ & {0.68} \\
      \bottomrule
    \end{tabular}
    \caption{Errors for one refinement in the solid domain.}
    \label{refinement_solid_1D}
  \end{center}
\end{table}

In Tables~\ref{uniform_1D}, \ref{refinement_fluid_1D} and \ref{refinement_solid_1D} we show the results for the one-dimensional model problem. We consider three cases, fully uniform time-stepping, a set-up with one refinement in time in the fluid domain only, and a similar configuration where only the time-steps in the solid domain are refined. For all of the simulations, the space mesh remains unchanged and only the time meshes are refined. We use a rather coarse space mesh with 100 cells only. To illustrate how the error is decomposed into different parts, we show the values of $\eta_{kh}^f$, $\eta_{kh}^s$ as well as $\eta_{kh}$. In Table~\ref{uniform_1D} we collect the numbers for fully uniform time-stepping. We can clearly see linear convergence of the error expected for the backward Euler time-stepping scheme. Only on higher refinement levels does the convergence rate slightly deteriorate due to the space mesh's coarseness. We can notice that the contributions from the solid problem dominate. In Table~\ref{refinement_fluid_1D} we once refine the time-steps in the fluid domain while leaving the time-steps in the solid domain unchanged. Unsurprisingly, this change does not lead to a decrease in the overall error. In Table~\ref{refinement_solid_1D} we instead refine the solid time steps which leads to a desirable result in a reduction of the overall error. In both cases of partial refinement, the error contributions from the unrefined subproblems slightly increase.

\subsubsection{Tensor-product space-time FEM results}
For the space-time finite element discretization, we use a Galerkin discretization with $\dG(1)$ in time and linear finite elements in space. 
For this numerical test, we use a penalty parameter $\gamma = 1000$, an initial coarse temporal mesh with $|\mathcal{T}_k^{\text{coarse}}| = 4$ temporal elements and 10 spatial degrees of freedom each for fluid and solid. 
To test the performance of temporal multirate space-time FEM, we uniformly refine the fluid or solid temporal mesh up to 8 times. This means that e.g. for each 1 solid temporal element we have 256 fluid temporal elements. The initial mesh is being refined uniformly in space and time. To measure the convergence, we compute the error between the analytical solution $U$ and the finite element solution $U_{kh}$ using the quantity of interest
\begin{align*}
    J(U - U_{kh}) &= \|U - U_{kh}\|_{L^2(I, L^2(\Omega))}.
\end{align*}

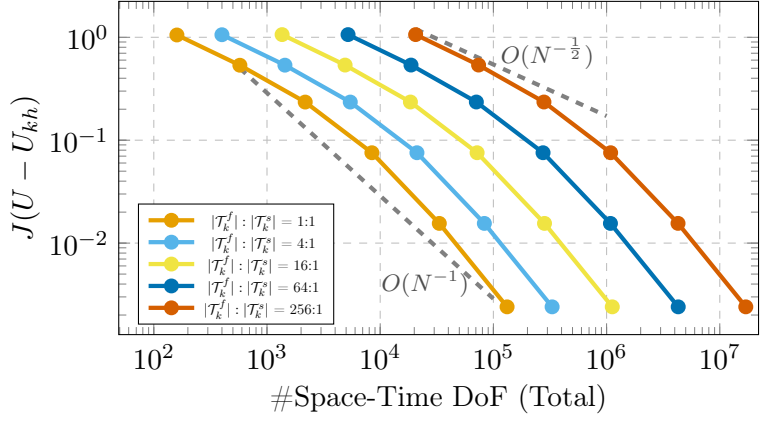
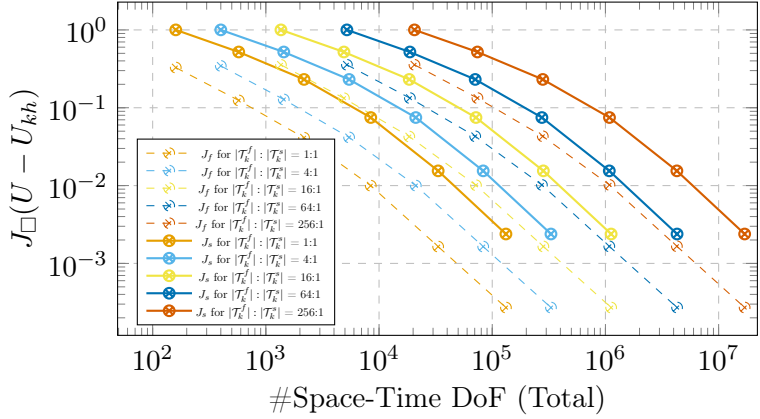
\begin{figure}[H]
    \centering
    \subfloat[Total error]{
        \centering
        \begin{tikzpicture}
            \begin{axis}[
                xlabel={\#Space-Time DoF {(Total)}},
                ylabel={$J(U - U_{kh})$},
                xmin=0.0, xmax=22000000,
                legend pos=south west,
                legend style={nodes={scale=0.5, transform shape}},
                xmajorgrids=true,
                ymajorgrids=true,
                grid style=dashed,
                xmode = log,
                ymode = log,
                width = 10cm, 
                height = 6cm,
            ]
        
            \addplot[
            color=orange,
            mark=otimes,
            style=ultra thick,
            ]
            coordinates {
            (160,1.053e+00)(576,5.351e-01)(2176,2.347e-01)(8448,7.554e-02)(33280,1.554e-02)(132096,2.402e-03)
            };
        
            \addplot[
            color=skyblue,
            mark=otimes,
            style=ultra thick,
            ]
            coordinates {
            (400,1.059e+00)(1440,5.359e-01)(5440,2.348e-01)(21120,7.555e-02)(83200,1.554e-02)(330240,2.402e-03)
            };
        
            \addplot[
            color=yellow,
            mark=otimes,
            style=ultra thick,
            ]
            coordinates {
            (1360,1.062e+00)(4896,5.371e-01)(18496,2.350e-01)(71808,7.556e-02)(282880,1.554e-02)(1122816,2.402e-03)
            };
        
            \addplot[
            color=blue,
            mark=otimes,
            style=ultra thick,
            ]
            coordinates {
            (5200,1.062e+00)(18720,5.372e-01)(70720,2.350e-01)(274560,7.557e-02)(1081600,1.554e-02)(4293120,2.402e-03)
            };
        
            \addplot[
            color=vermillion,
            mark=otimes,
            style=ultra thick,
            ]
            coordinates {
            (20560,1.062e+00)(74016,5.372e-01)(279616,2.350e-01)(1085568,7.557e-02)(4276480,1.555e-02)(16974336,2.402e-03)
            }; 
        
            \addplot[
            color=gray,
            style=ultra thick,
            no markers,
            dashed
            ]
            coordinates {
            (20560, 1.2)(1000000,0.17206510395777522) 
            };
            \node[] at (axis cs: 300000, 0.7) {\color{gray!60!black}\footnotesize $O(N^{-\frac{1}{2}})$};
        
            \addplot[
            color=gray,
            style=ultra thick,
            no markers,
            dashed
            ]
            coordinates {
            (576, 0.5)(100000,0.00288) 
            };
            \node[] at (axis cs: 25000, 0.004) {\color{gray!60!black}\footnotesize $O(N^{-1})$};
            
            \legend{{} {$|\mathcal{T}_k^f| : |\mathcal{T}_k^s|$ = 1:1},{} {$|\mathcal{T}_k^f| : |\mathcal{T}_k^s|$ = 4:1},{} {$|\mathcal{T}_k^f| : |\mathcal{T}_k^s|$ = 16:1},{} {$|\mathcal{T}_k^f| : |\mathcal{T}_k^s|$ = 64:1},{} {$|\mathcal{T}_k^f| : |\mathcal{T}_k^s|$ = 256:1},}
            \end{axis}
    \end{tikzpicture}
    }
    \hfill
    \subfloat[\auth{Fluid and solid errors}]{
        \centering
        \begin{tikzpicture}
            \begin{axis}[
                xlabel={\#Space-Time DoF {(Total)}},
                ylabel={$J_\Box(U - U_{kh})$},
                xmin=0.0, xmax=22000000,
                legend pos=south west,
                legend style={nodes={scale=0.4, transform shape}},
                xmajorgrids=true,
                ymajorgrids=true,
                grid style=dashed,
                xmode = log,
                ymode = log,
                width = 10cm, 
                height = 6cm,
            ]
            \addplot[
                color=orange,
            mark=otimes,
                style=dashed,
                ]
                coordinates {
            (160,3.282e-01)(576,1.242e-01)(2176,4.116e-02)(8448,9.967e-03)(33280,1.644e-03)(132096,2.708e-04)
            };
        
                \addplot[
                color=skyblue,
            mark=otimes,
                style=dashed,
                ]
                coordinates {
            (400,3.446e-01)(1440,1.278e-01)(5440,4.159e-02)(21120,9.996e-03)(83200,1.645e-03)(330240,2.708e-04)
            };
        
                \addplot[
                color=yellow,
            mark=otimes,
                style=dashed,
                ]
                coordinates {
            (1360,3.556e-01)(4896,1.325e-01)(18496,4.263e-02)(71808,1.010e-02)(282880,1.650e-03)(1122816,2.709e-04)
            };
        
                \addplot[
                color=blue,
            mark=otimes,
                style=dashed,
                ]
                coordinates {
            (5200,3.561e-01)(18720,1.331e-01)(70720,4.294e-02)(274560,1.019e-02)(1081600,1.657e-03)(4293120,2.711e-04)
            };
        
                \addplot[
                color=vermillion,
            mark=otimes,
                style=dashed,
                ]
                coordinates {
            (20560,3.561e-01)(74016,1.331e-01)(279616,4.296e-02)(1085568,1.020e-02)(4276480,1.660e-03)(16974336,2.713e-04)
            };
        
                \addplot[
                color=orange,
            mark=otimes,
                style=thick,
                ]
                coordinates {
            (160,1.001e+00)(576,5.205e-01)(2176,2.311e-01)(8448,7.488e-02)(33280,1.546e-02)(132096,2.387e-03)
            };
        
                \addplot[
                color=skyblue,
            mark=otimes,
                style=thick,
                ]
                coordinates {
            (400,1.001e+00)(1440,5.205e-01)(5440,2.311e-01)(21120,7.488e-02)(83200,1.546e-02)(330240,2.387e-03)
            };
        
                \addplot[
                color=yellow,
            mark=otimes,
                style=thick,
                ]
                coordinates {
            (1360,1.001e+00)(4896,5.205e-01)(18496,2.311e-01)(71808,7.488e-02)(282880,1.546e-02)(1122816,2.387e-03)
            };
        
                \addplot[
                color=blue,
            mark=otimes,
                style=thick,
                ]
                coordinates {
            (5200,1.001e+00)(18720,5.205e-01)(70720,2.311e-01)(274560,7.488e-02)(1081600,1.546e-02)(4293120,2.387e-03)
            };
        
                \addplot[
                color=vermillion,
            mark=otimes,
                style=thick,
                ]
                coordinates {
            (20560,1.001e+00)(74016,5.205e-01)(279616,2.311e-01)(1085568,7.488e-02)(4276480,1.546e-02)(16974336,2.387e-03)
            };

            \legend{{$J_f$ for $|\mathcal{T}_k^f| : |\mathcal{T}_k^s|$ = 1:1}, {$J_f$ for $|\mathcal{T}_k^f| : |\mathcal{T}_k^s|$ = 4:1}, {$J_f$ for $|\mathcal{T}_k^f| : |\mathcal{T}_k^s|$ = 16:1}, {$J_f$ for $|\mathcal{T}_k^f| : |\mathcal{T}_k^s|$ = 64:1}, {$J_f$ for $|\mathcal{T}_k^f| : |\mathcal{T}_k^s|$ = 256:1}, {$J_s$ for $|\mathcal{T}_k^f| : |\mathcal{T}_k^s|$ = 1:1}, {$J_s$ for $|\mathcal{T}_k^f| : |\mathcal{T}_k^s|$ = 4:1}, {$J_s$ for $|\mathcal{T}_k^f| : |\mathcal{T}_k^s|$ = 16:1}, {$J_s$ for $|\mathcal{T}_k^f| : |\mathcal{T}_k^s|$ = 64:1}, {$J_s$ for $|\mathcal{T}_k^f| : |\mathcal{T}_k^s|$ = 256:1},}
            \end{axis}
    \end{tikzpicture}}
    \caption{Convergence plots for the 1+1D heat and wave equation problem with a finer fluid temporal mesh and uniform refinement in space and time}
    \label{fig:convergence_space_time_1d_finer_fluid}
\end{figure}

Using a finer fluid temporal mesh, we get the convergence plot in Figure \ref{fig:convergence_space_time_1d_finer_fluid}. Therein, we observe that using a finer temporal mesh for the fluid in comparison to the solid does not improve the error between the finite element and the analytical solution.
\auth{Moreover, in Figure \ref{fig:convergence_space_time_1d_finer_fluid}(b) we observe that the error in the solid domain $J_s$ (thick lines) is almost an order of magnitude larger than the error in the fluid domain $J_f$ (dashed lines).}

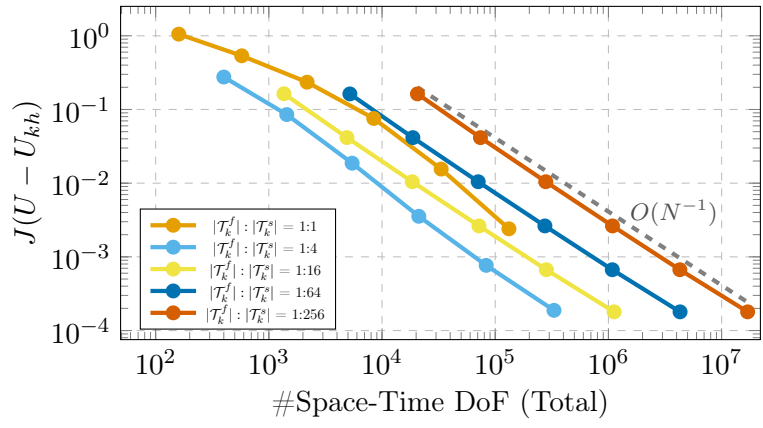
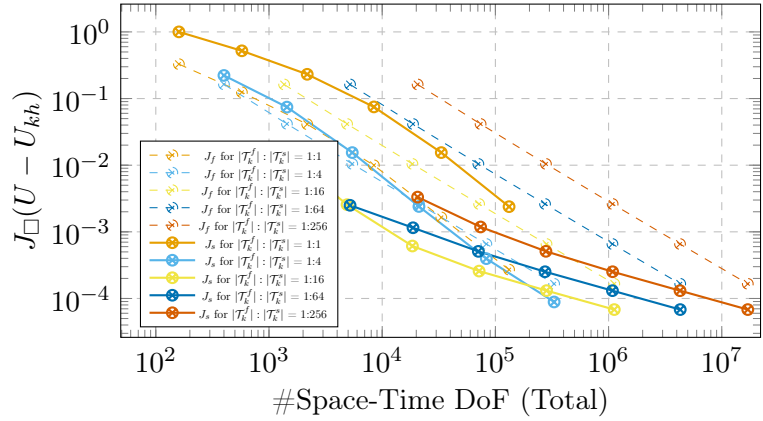
\begin{figure}[H]%
    \centering
    \subfloat[Total error]{
        \begin{tikzpicture}
        \begin{axis}[
            xlabel={\#Space-Time DoF {(Total)}},
            ylabel={$J(U - U_{kh})$},
            xmin=0.0, xmax=22000000,
            legend pos=south west,
            legend style={nodes={scale=0.5, transform shape}},
            xmajorgrids=true,
            ymajorgrids=true,
            grid style=dashed,
            xmode = log,
            ymode = log,
            width = 10cm, 
            height = 6cm,
        ]

        \addplot[
        color=orange,
        mark=otimes,
        style=ultra thick,
        ]
        coordinates {
        (160,1.053e+00)(576,5.351e-01)(2176,2.347e-01)(8448,7.554e-02)(33280,1.554e-02)(132096,2.402e-03)
        };

        \addplot[
        color=skyblue,
        mark=otimes,
        style=ultra thick,
        ]
        coordinates {
        (400,2.757e-01)(1440,8.531e-02)(5440,1.867e-02)(21120,3.552e-03)(83200,7.667e-04)(330240,1.881e-04)
        };

        \addplot[
        color=yellow,
        mark=otimes,
        style=ultra thick,
        ]
        coordinates {
        (1360,1.635e-01)(4896,4.157e-02)(18496,1.046e-02)(71808,2.629e-03)(282880,6.690e-04)(1122816,1.794e-04)
        };

        \addplot[
        color=blue,
        mark=otimes,
        style=ultra thick,
        ]
        coordinates {
        (5200,1.627e-01)(18720,4.151e-02)(70720,1.046e-02)(274560,2.629e-03)(1081600,6.689e-04)(4293120,1.794e-04)
        };

        \addplot[
        color=vermillion,
        mark=otimes,
        style=ultra thick,
        ]
        coordinates {
        (20560,1.628e-01)(74016,4.151e-02)(279616,1.046e-02)(1085568,2.629e-03)(4276480,6.689e-04)(16974336,1.794e-04)
        };

        \addplot[
        color=gray,
        style=ultra thick,
        no markers,
        dashed
        ]
        coordinates {
        (20560, 0.2)(18000000,0.00022844444444444445) 
        };
        \node[] at (axis cs: 3800000, 0.004) {\color{gray!60!black}\footnotesize $O(N^{-1})$};
        
        \legend{{} {$|\mathcal{T}_k^f| : |\mathcal{T}_k^s|$ = 1:1},{} {$|\mathcal{T}_k^f| : |\mathcal{T}_k^s|$ = 1:4},{} {$|\mathcal{T}_k^f| : |\mathcal{T}_k^s|$ = 1:16},{} {$|\mathcal{T}_k^f| : |\mathcal{T}_k^s|$ = 1:64},{} {$|\mathcal{T}_k^f| : |\mathcal{T}_k^s|$ = 1:256},}
        \end{axis}
        \end{tikzpicture}}
    \hfill
    \subfloat[\auth{Fluid and solid errors}]{
        \begin{tikzpicture}
            \begin{axis}[
                xlabel={\#Space-Time DoF {(Total)}},
                ylabel={$J_\Box(U - U_{kh})$},
                xmin=0.0, xmax=22000000,
                legend pos=south west,
                legend style={nodes={scale=0.4, transform shape}},
                xmajorgrids=true,
                ymajorgrids=true,
                grid style=dashed,
                xmode = log,
                ymode = log,
                width = 10cm, 
                height = 6cm,
            ]
    
            \addplot[
                color=orange,
            mark=otimes,
                style=dashed,
                ]
                coordinates {
            (160,3.282e-01)(576,1.242e-01)(2176,4.116e-02)(8448,9.967e-03)(33280,1.644e-03)(132096,2.708e-04)
            };
        
                \addplot[
                color=skyblue,
            mark=otimes,
                style=dashed,
                ]
                coordinates {
            (400,1.634e-01)(1440,4.169e-02)(5440,1.046e-02)(21120,2.618e-03)(83200,6.564e-04)(330240,1.659e-04)
            };
        
                \addplot[
                color=yellow,
            mark=otimes,
                style=dashed,
                ]
                coordinates {
            (1360,1.627e-01)(4896,4.149e-02)(18496,1.044e-02)(71808,2.616e-03)(282880,6.560e-04)(1122816,1.659e-04)
            };
        
                \addplot[
                color=blue,
            mark=otimes,
                style=dashed,
                ]
                coordinates {
            (5200,1.627e-01)(18720,4.149e-02)(70720,1.044e-02)(274560,2.616e-03)(1081600,6.560e-04)(4293120,1.659e-04)
            };
        
                \addplot[
                color=vermillion,
            mark=otimes,
                style=dashed,
                ]
                coordinates {
            (20560,1.627e-01)(74016,4.149e-02)(279616,1.044e-02)(1085568,2.616e-03)(4276480,6.560e-04)(16974336,1.659e-04)
            };
        
                \addplot[
                color=orange,
            mark=otimes,
                style=thick,
                ]
                coordinates {
            (160,1.001e+00)(576,5.205e-01)(2176,2.311e-01)(8448,7.488e-02)(33280,1.546e-02)(132096,2.387e-03)
            };
        
                \addplot[
                color=skyblue,
            mark=otimes,
                style=thick,
                ]
                coordinates {
            (400,2.221e-01)(1440,7.443e-02)(5440,1.546e-02)(21120,2.401e-03)(83200,3.962e-04)(330240,8.863e-05)
            };
        
                \addplot[
                color=yellow,
            mark=otimes,
                style=thick,
                ]
                coordinates {
            (1360,1.539e-02)(4896,2.558e-03)(18496,6.127e-04)(71808,2.577e-04)(282880,1.312e-04)(1122816,6.816e-05)
            };
        
                \addplot[
                color=blue,
            mark=otimes,
                style=thick,
                ]
                coordinates {
            (5200,2.502e-03)(18720,1.151e-03)(70720,5.092e-04)(274560,2.520e-04)(1081600,1.310e-04)(4293120,6.815e-05)
            };
        
                \addplot[
                color=vermillion,
            mark=otimes,
                style=thick,
                ]
                coordinates {
            (20560,3.321e-03)(74016,1.183e-03)(279616,5.098e-04)(1085568,2.520e-04)(4276480,1.310e-04)(16974336,6.815e-05)
            };
                                
            \legend{{$J_f$ for $|\mathcal{T}_k^f| : |\mathcal{T}_k^s|$ = 1:1}, {$J_f$ for $|\mathcal{T}_k^f| : |\mathcal{T}_k^s|$ = 1:4}, {$J_f$ for $|\mathcal{T}_k^f| : |\mathcal{T}_k^s|$ = 1:16}, {$J_f$ for $|\mathcal{T}_k^f| : |\mathcal{T}_k^s|$ = 1:64}, {$J_f$ for $|\mathcal{T}_k^f| : |\mathcal{T}_k^s|$ = 1:256}, {$J_s$ for $|\mathcal{T}_k^f| : |\mathcal{T}_k^s|$ = 1:1}, {$J_s$ for $|\mathcal{T}_k^f| : |\mathcal{T}_k^s|$ = 1:4}, {$J_s$ for $|\mathcal{T}_k^f| : |\mathcal{T}_k^s|$ = 1:16}, {$J_s$ for $|\mathcal{T}_k^f| : |\mathcal{T}_k^s|$ = 1:64}, {$J_s$ for $|\mathcal{T}_k^f| : |\mathcal{T}_k^s|$ = 1:256},}
            \end{axis}
            \end{tikzpicture}
    }
    \caption{Convergence plot for the 1+1D heat and wave equation problem with a finer solid temporal mesh and uniform refinement in space and time}
    \label{fig:convergence_space_time_1d_finer_solid}
\end{figure}

Repeating the same convergence tests again for a finer solid temporal mesh, the results are shown in Figure \ref{fig:convergence_space_time_1d_finer_solid}.
Here, we observe that using a finer temporal mesh for the solid in comparison to the fluid greatly reduces the error between the finite element and the analytical solution. 
In particular, comparing 1:1 and 1:4 temporal meshes, we observe that the error on the coarsest grid is four times lower when the temporal mesh of the solid is finer and we then get also optimal convergence rates.
Using even finer temporal meshes, e.g. 1:16, further reduces the error but going beyond that, i.e. 1:64 and 1:256, almost does not improve the finite element solution anymore. 
\auth{Another explanation of these observations is given in Figure \ref{fig:convergence_space_time_1d_finer_solid}(b), where we decompose the error into its fluid and solid parts.
For the 1:1 temporal mesh, the error in the solid domain $J_s$ (thick lines) is almost an order of magnitude larger than the error in the fluid domain $J_f$ (dashed lines).
For the 1:4 temporal mesh, both lines are close to each other, which signifies that the error in the fluid and solid domain are of the same order of magnitude.
Then for the 1:16 temporal mesh and finer meshes, the error in the fluid domain is between one and two orders of magnitude larger than the error in the solid domain.
Hence, for this problem configuration it seems optimal to use four times as many solid temporal elements as fluid temporal elements. 
}

\begin{figure}[H]
     \centering
     \subfloat[Displacement]{
    \includegraphics[width=4cm]{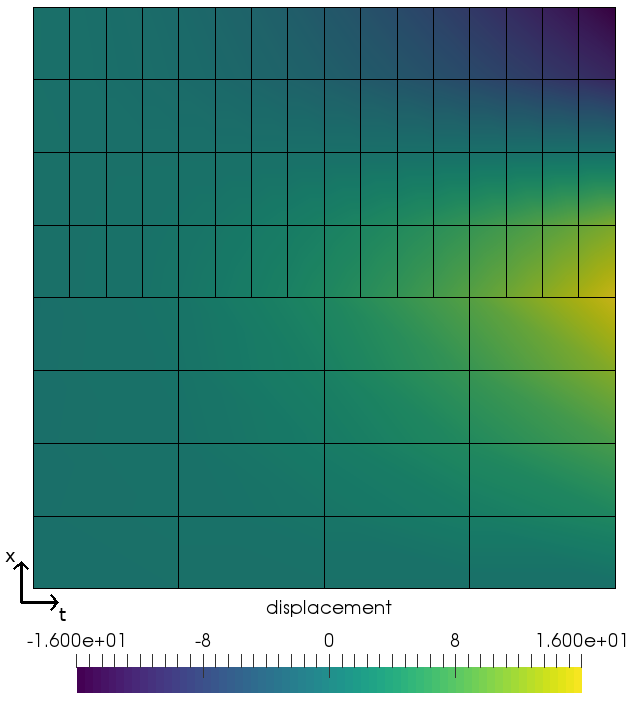}}%
    \hspace{2cm}
     \subfloat[Velocity]{
    \includegraphics[width=4cm]{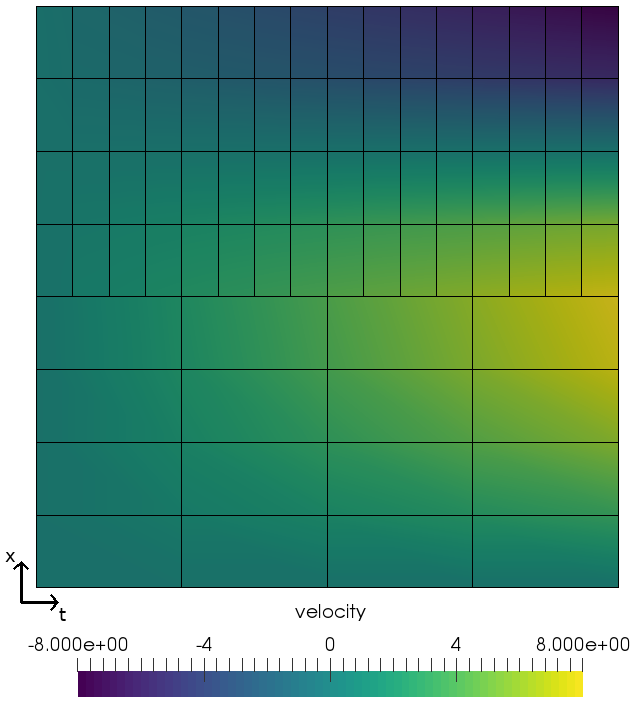}}%
    \\
    \subfloat[Displacement error]{
    \includegraphics[width=4cm]{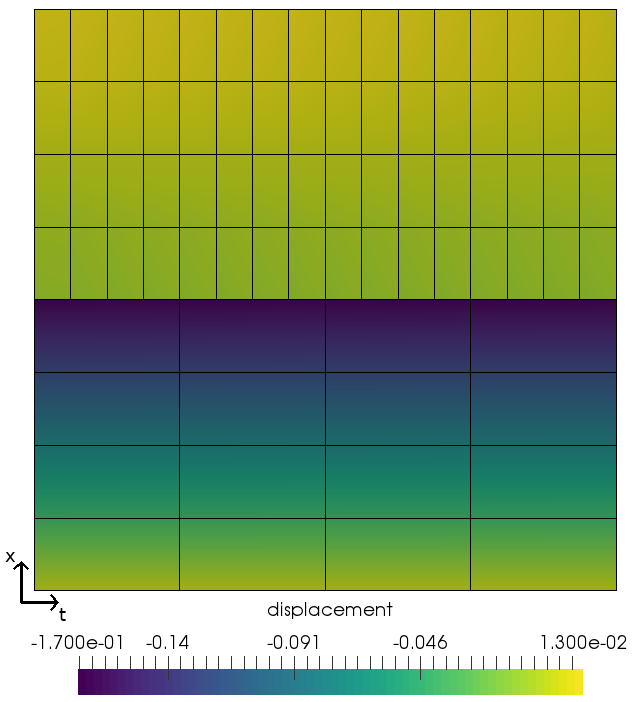}}%
    \hspace{2cm}
     \subfloat[Velocity error]{
    \includegraphics[width=4cm]{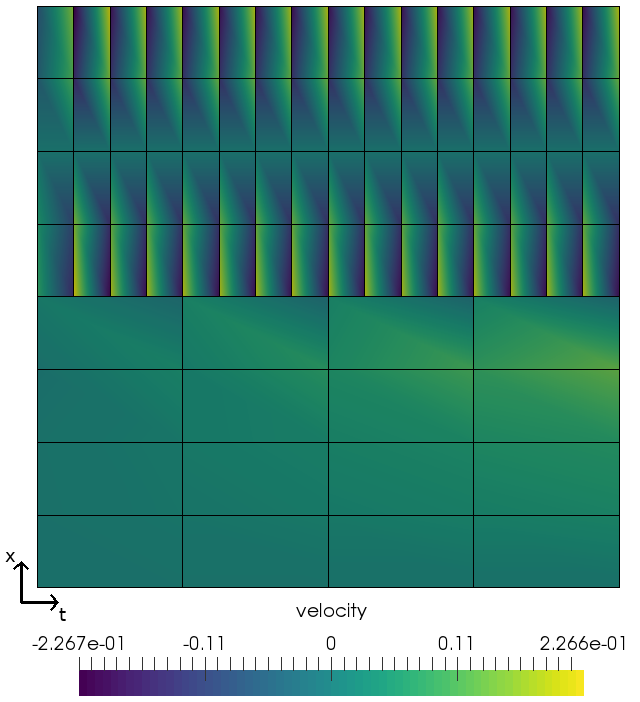}}%
    \caption{Space-time FEM solution for the 1+1D heat and wave equation problem with 4 solid temporal elements per 1 fluid temporal element}
    \label{fig:space_time_solution_1d}
\end{figure}

In Figure \ref{fig:space_time_solution_1d}, we show an example of a space-time FEM solution on the coarsest spatial and temporal meshes. Here, the temporal mesh of the solid, shown in the upper half of the space-time domain, is four times as fine as for the fluid, which is plotted in the lower half of the space-time domain. Additionally, we visualize the error to the analytical solution and we observe especially in the velocity error that a discontinuous time discretization has been used.

\begin{figure}[H]
     \centering
     \subfloat[Displacement]{
    \includegraphics[width=8cm]{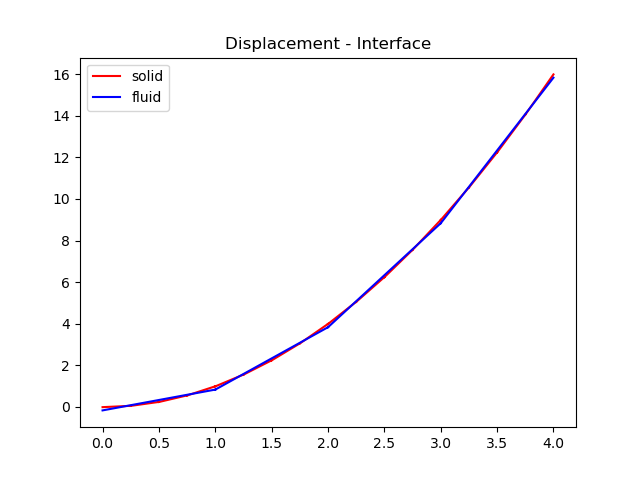}}%
     \subfloat[Velocity]{
    \includegraphics[width=8cm]{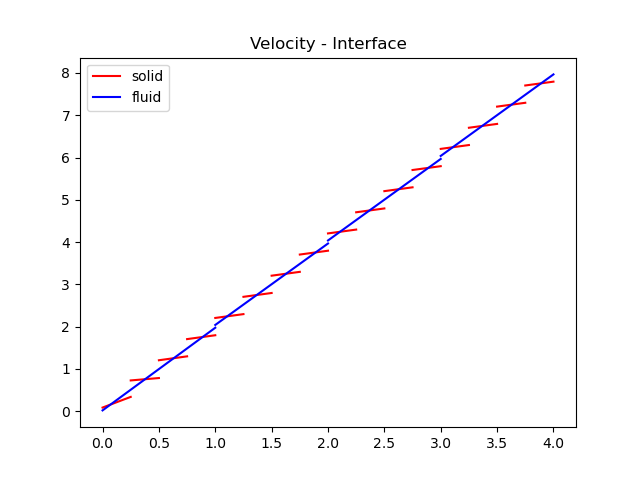}}%
    \caption{Temporal evolution at the interface of the space-time FEM solution for the 1+1D heat and wave equation problem with 4 solid temporal elements per 1 fluid temporal element}
    \label{fig:space_time_interface_1d}
\end{figure}

In Figure \ref{fig:space_time_interface_1d}, we now visualize the temporal evolution of the continuity of the displacement and the velocity at the interface of the space-time finite element solution from Figure \ref{fig:space_time_solution_1d}. 
Due to the different number of temporal elements for fluid and solid, we cannot expect a perfect match between the solutions at the interface. 
Nevertheless, the solid and the fluid displacement almost coincide. The solid and the fluid velocity are not too different either, but here we again observe the jumps in the solid solution between the temporal elements.

\subsection{2+1D heat and wave equation \auth{with fluid source term}}
\label{sec_2d_heat_wave}

As the second numerical test, we again consider a coupled heat and wave equation from Section \ref{sec_model_problem_heatwave} but now in two spatial dimensions. For this we use Configuration 2.1 from \cite{Soszynska2021}. The time interval is $I:=(0,1)$, the spatial domain for the heat equation - the fluid domain - is $\Omega_f := (0,4) \times (0,1)$ and the spatial domain for the wave equation - the solid domain - is $\Omega_s := (0,4) \times (-1,0)$. The interface is given by $\Gamma := \bar{\Omega}_f \cap \bar{\Omega}_s = (0,4) \times \{0\}.$
For our test we choose the parameters in the PDEs to be 
\begin{align*}
    \nu = 0.001, \quad \beta = \begin{pmatrix}
    2 \\ 0
\end{pmatrix}, \quad \lambda = 1000,\quad \delta = 0.1.   
\end{align*}
As right hand sides we choose
\begin{align*}
    g_s(x, t) &= 0, \\
    g_f(x, t) &= \begin{cases}
        e^{-\left((x_1-\frac{1}{2})^2+(x_2-\frac{1}{2})^2\right)}  & \text{if } 0 \leq t \leq 0.1,\\
        0 & \text{else},
    \end{cases}
\end{align*}
\auth{i.e. we only have a source term in the fluid domain.}
The boundaries are defined as $\Gamma_D^f = (0,4) \times \{1\}$, $\Gamma_N^f = \{0,4\} \times (0,1)$, $\Gamma_D^s = \{0,4\} \times (-1,0)$, $\Gamma_N^s =(0,4) \times \{-1\}$ and the variational formulation is shown in (\ref{eq:variational_form_heat_wave}).

\auth{From this subsection on, we only employ the} space-time finite element discretization (but not anymore the time marching schemes) in which 
we use a Galerkin discretization with $\dG(1)$ in time and linear finite elements in space. 
For this numerical test, we use as penalty parameter $\gamma = 1000$, a spatial mesh with 80 spatial cells in the x-direction and 20 spatial cells in the y-direction, and an initial coarse temporal mesh with $|\mathcal{T}_k^{\text{coarse}}| = 50$ temporal elements.
To test the performance of temporal multirate space-time FEM, we uniformly refine the fluid or solid temporal mesh up to 4 times. This means that e.g. for each 1 solid temporal element we have 16 fluid temporal elements. The initial mesh is being refined uniformly only in time. To measure the convergence, we consider the quantity of interest 
\begin{equation*}
    \auth{J(U_{kh}) = \nu \|\nabla_x v_{kh}\|^2_{L^2(I, L^2(\Omega_f))},}
\end{equation*}
and compare with the reference value $J(U) := 2.48587692 \cdot 10^{-4}$ from a simulation with $50,000$ temporal elements.

\begin{figure}[H]%
    \centering
    \begin{tikzpicture}
    \begin{axis}[
        xlabel={$|\mathcal{T}_k^{\text{coarse}}|$},
        ylabel={$J(U) - J(U_{kh})$},
        xmin=0.0, xmax=500,
        legend pos= north east, 
        legend style={nodes={scale=1, transform shape}},
        xmajorgrids=true,
        ymajorgrids=true,
        grid style=dashed,
        xmode = log,
        ymode = log,
        width = 15cm, 
        height = 9cm,
        xtick = {50,100,200,400},
        xticklabels = {50,100,200,400},
    ]

    \addplot[
    color=orange,
	mark=otimes,
    style=ultra thick,
    ]
    coordinates {
	(50,1.702211997061959e-09)(100,2.579022337486844e-10)(200,3.943867649992827e-11)(400,5.522828722126505e-12)
	};

    \addplot[
    color=skyblue,
	mark=otimes,
    style=ultra thick,
    ]
    coordinates {
	(50,2.637860456325761e-10)(100,5.517309059708403e-11)(200,9.688701663873589e-12)(400,1.352230718891512e-12)
	};

    \addplot[
    color=yellow,
	mark=otimes,
    style=ultra thick,
    ]
    coordinates {
	(50,6.128633527318489e-11)(100,2.545320546129201e-11)(200,5.518989345393299e-12)(400,7.973878718772753e-13)
	};

    \addplot[
    color=blue,
	mark=otimes,
    style=ultra thick,
    ]
    coordinates {
	(50,3.158910432494652e-11)(100,2.128758980072046e-11)(200,4.964151539919165e-12)(400,7.255025573360552e-13)
	};

    \addplot[
    color=vermillion,
	mark=otimes,
    style=ultra thick,
    ]
    coordinates {
	(50,2.742656080544081e-11)(100,2.073323457363330e-11)(200,4.892428096762297e-12)(400,7.146690465982541e-13)
	};

    \addplot[
    color=gray,
    style=ultra thick,
    no markers,
    dashed
    ]
    coordinates {
    (50,4e-09)(400,7.8125e-12) 
    };
    \node[] at (axis cs: 300, 4e-11) {\color{gray!60!black}\footnotesize $O(N^{-3})$};

    \legend{{} {$|\mathcal{T}_k^f| : |\mathcal{T}_k^s|$ = 1:1},{} {$|\mathcal{T}_k^f| : |\mathcal{T}_k^s|$ = 2:1},{} {$|\mathcal{T}_k^f| : |\mathcal{T}_k^s|$ = 4:1},{} {$|\mathcal{T}_k^f| : |\mathcal{T}_k^s|$ = 8:1},{} {$|\mathcal{T}_k^f| : |\mathcal{T}_k^s|$ = 16:1},}
    \end{axis}
    \end{tikzpicture}
    \caption{Convergence plot for the 2+1D heat and wave equation problem with a finer fluid temporal mesh and uniform refinement in time}
    \label{fig:convergence_space_time_2d_finer_fluid}
\end{figure}
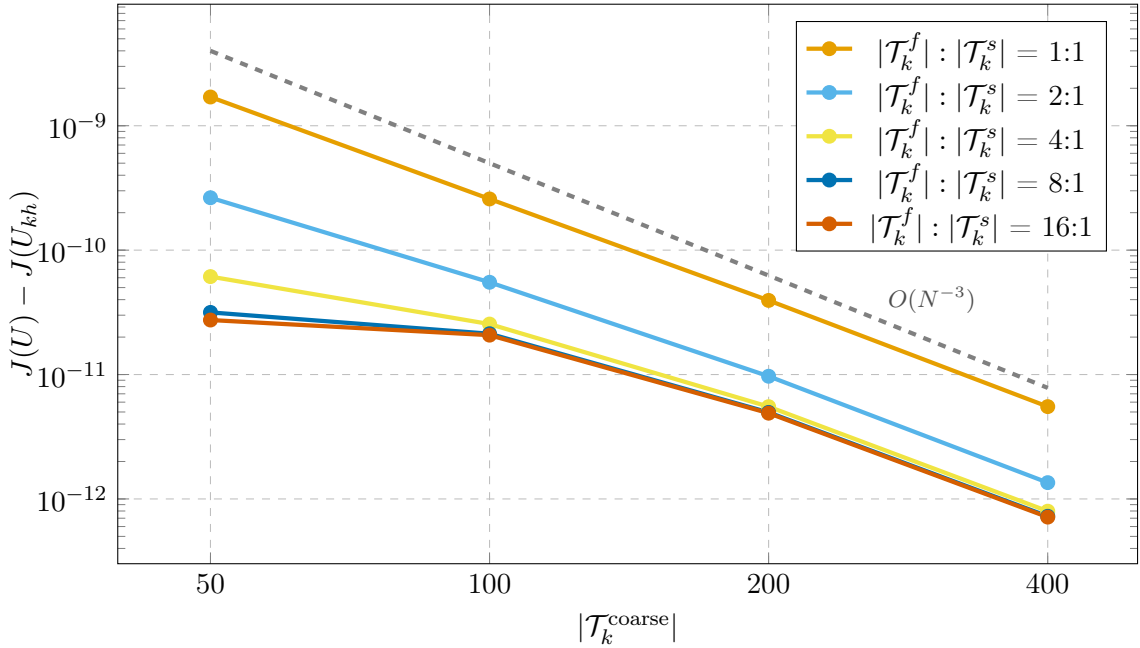

Using a finer fluid temporal mesh, we get the convergence plot in Figure \ref{fig:convergence_space_time_2d_finer_fluid}. Therein, we observe that using a finer temporal mesh for the fluid in comparison to the solid greatly reduces the error between the finite element quantity of interest and its reference value. These improvements are greatest for 1:1, 2:1 and 4:1 temporal meshes. Even finer fluid temporal meshes have only marginal benefit, since although they initially reduce the error in the goal functional, on finer temporal meshes there are minor differences between the 4:1 and 16:1 temporal meshes.

\begin{figure}[H]%
    \centering
    \begin{tikzpicture}
    \begin{axis}[
        xlabel={$|\mathcal{T}_k^{\text{coarse}}|$},
        ylabel={$J(U) - J(U_{kh})$},
        xmin=0.0, xmax=500,
        legend pos= north east, 
        legend style={nodes={scale=1, transform shape}},
        xmajorgrids=true,
        ymajorgrids=true,
        grid style=dashed,
        xmode = log,
        ymode = log,
        width = 15cm, 
        height = 9cm,
        xtick = {50,100,200,400},
        xticklabels = {50,100,200,400},
    ]

    \addplot[
    color=orange,
	mark=otimes,
    style=ultra thick,
    ]
    coordinates {
	(50,1.702211997061959e-09)(100,2.579022337486844e-10)(200,3.943867649992827e-11)(400,5.522828722126505e-12)
	};

    \addplot[
    color=skyblue,
	mark=otimes,
    style=ultra thick,
    ]
    coordinates {
	(50,1.697079597294190e-09)(100,2.423024584251830e-10)(200,3.527609286494218e-11)(400,4.893984523191008e-12)
	};

    \addplot[
    color=yellow,
	mark=otimes,
    style=ultra thick,
    ]
    coordinates {
	(50,1.682290524116295e-09)(100,2.381803292463881e-10)(200,3.464679519844804e-11)(400,4.812149539522997e-12)
	};

    \addplot[
    color=blue,
	mark=otimes,
    style=ultra thick,
    ]
    coordinates {
	(50,1.678501410236864e-09)(100,2.375603154946768e-10)(200,3.456478912767721e-11)(400,4.801838288971683e-12)
	};

    \addplot[
    color=vermillion,
	mark=otimes,
    style=ultra thick,
    ]
    coordinates {
	(50,1.677957375963938e-09)(100,2.374795938782197e-10)(200,3.455440870502730e-11)(400,4.800559960400214e-12)
	};

    \addplot[
    color=gray,
    style=ultra thick,
    no markers,
    dashed
    ]
    coordinates {
    (50,4e-09)(400,7.8125e-12) 
    };
    \node[] at (axis cs: 300, 4e-11) {\color{gray!60!black}\footnotesize $O(N^{-3})$};
 
    \legend{{} {$|\mathcal{T}_k^f| : |\mathcal{T}_k^s|$ = 1:1},{} {$|\mathcal{T}_k^f| : |\mathcal{T}_k^s|$ = 1:2},{} {$|\mathcal{T}_k^f| : |\mathcal{T}_k^s|$ = 1:4},{} {$|\mathcal{T}_k^f| : |\mathcal{T}_k^s|$ = 1:8},{} {$|\mathcal{T}_k^f| : |\mathcal{T}_k^s|$ = 1:16},}

    \end{axis}
    \end{tikzpicture}
    \caption{Convergence plot for the 2+1D heat and wave equation problem with a finer solid temporal mesh and uniform refinement in time}
    \label{fig:convergence_space_time_2d_finer_solid}
\end{figure}
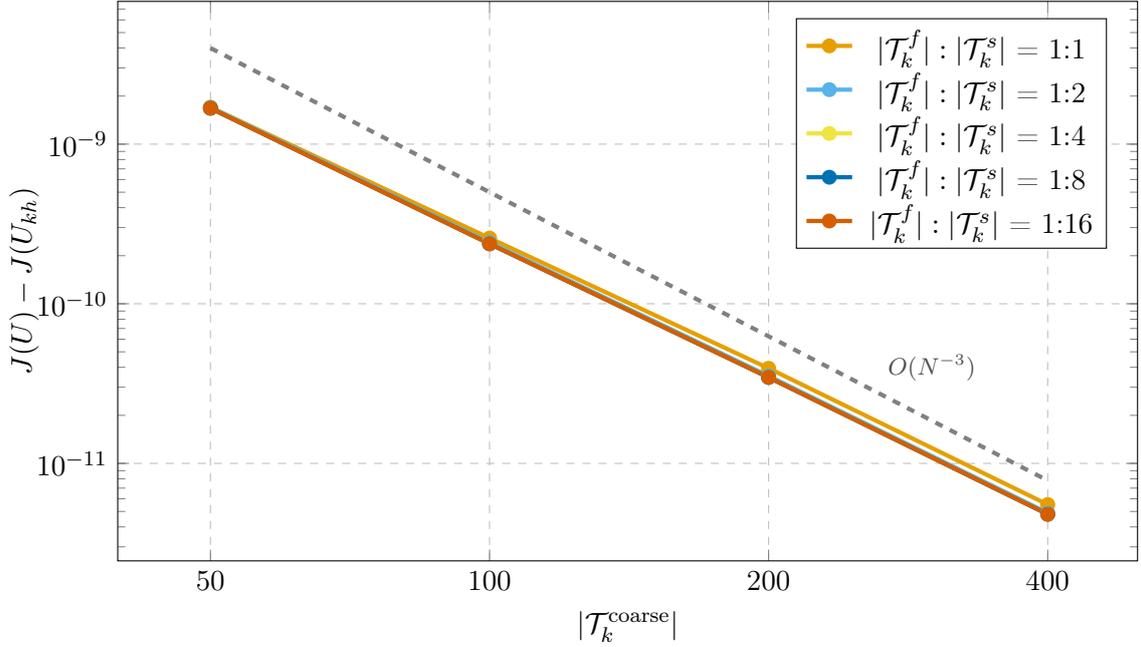

Performing the same convergence test for a finer solid temporal mesh, the results are shown in Figure \ref{fig:convergence_space_time_2d_finer_solid}. Here, we observe that using a finer temporal mesh for the solid in comparison to the fluid does not improve the error in the quantity of interest. This was to be expected since we only have a non-zero right hand side in the fluid domain and thus the solution behavior is mainly dictated by the heat equation. 

\subsection{2+1D heat and wave equation with solid source term}
\label{sec_2d_heat_wave_solid_source}
{
As the third numerical test, we again consider a two-dimensional coupled heat and wave equation, but here we use Configuration 2.2 from \cite{Soszynska2021}. The computational space-time domain and the physical parameters are the same as in Section \ref{sec_2d_heat_wave}. The only difference is that we now have a solid source term in the wave equation. More concretely, as right hand sides we now choose 
\begin{align*}
    g_f(x, t) &= 0, \\
    g_s(x, t) &= \begin{cases}
        e^{-\left((x_1-\frac{1}{2})^2+(x_2+\frac{1}{2})^2\right)}  & \text{if } 0 \leq t \leq 0.1,\\
        0 & \text{else}.
    \end{cases}
\end{align*}
We employ a space-time finite element discretization with $\dG(1)$ in time and linear finite elements in space for this numerical test.
The simulations are performed on a coarse space mesh consisting of 80 spatial cells in the x-direction and 20 spatial cells in the y-direction, and an initial coarse temporal mesh with $|\mathcal{T}_k^{\text{coarse}}| = 50$ temporal elements.
To test the performance of temporal multirate space-time FEM, we uniformly refine the fluid or solid temporal mesh up to 4 times. To measure the convergence, we consider the quantity of interest 
\begin{align*}
    J(U_{kh}) = \lambda \|\nabla_x u_{kh}\|^2_{L^2(I, L^2(\Omega_s))},
\end{align*}
and compare with the reference value $J(U) :=  7.14276824 \cdot 10^{-4}$ which is taken from computations on a fine temporal mesh with $50,000$ temporal elements.}

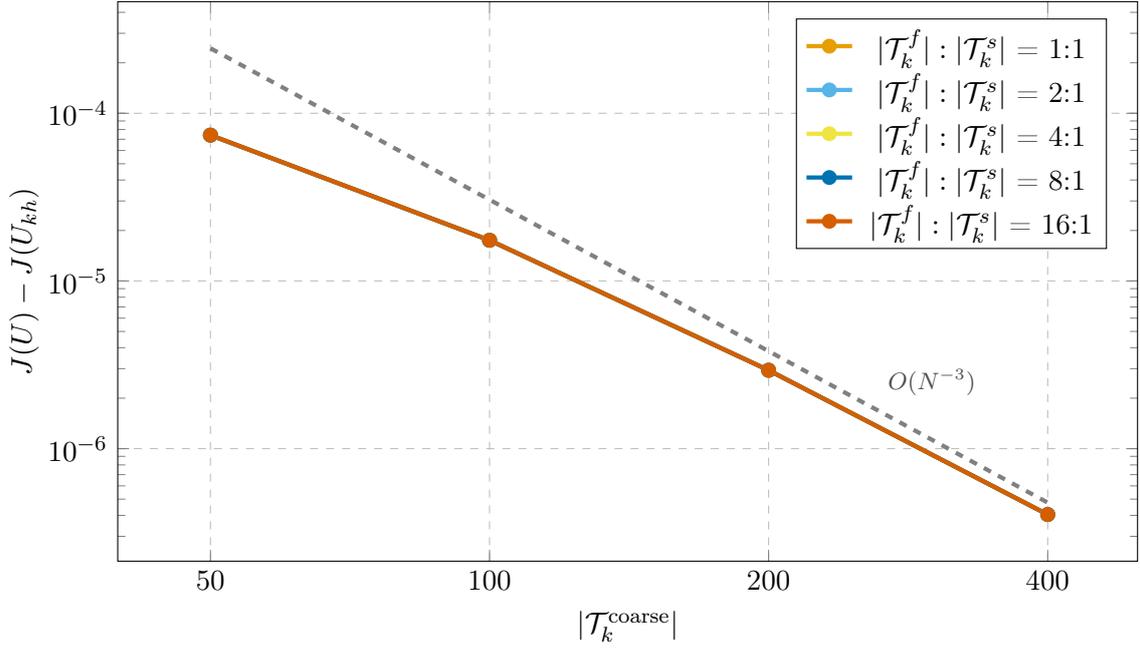
\begin{figure}[H]%
    \centering
    \begin{tikzpicture}
    \begin{axis}[
        xlabel={$|\mathcal{T}_k^{\text{coarse}}|$},
        ylabel={$J(U) - J(U_{kh})$},
        xmin=0.0, xmax=500,
        legend pos= north east, 
        legend style={nodes={scale=1, transform shape}},
        xmajorgrids=true,
        ymajorgrids=true,
        grid style=dashed,
        xmode = log,
        ymode = log,
        width = 15cm, 
        height = 9cm,
        xtick = {50,100,200,400},
        xticklabels = {50,100,200,400},
    ]

    \addplot[
        color=orange,
        mark=otimes,
        style=ultra thick,
        ]
        coordinates {
        (50,7.407397938424647e-05)(100,1.747727249566005e-05)(200,2.934075263809842e-06)(400,4.046231868823792e-07)
        };
    
        \addplot[
        color=skyblue,
        mark=otimes,
        style=ultra thick,
        ]
        coordinates {
        (50,7.407462724287465e-05)(100,1.747749454746505e-05)(200,2.934118730785087e-06)(400,4.046291931052420e-07)
        };
    
        \addplot[
        color=yellow,
        mark=otimes,
        style=ultra thick,
        ]
        coordinates {
        (50,7.407464420046558e-05)(100,1.747751859110963e-05)(200,2.934124546065202e-06)(400,4.046299988749724e-07)
        };
    
        \addplot[
        color=blue,
        mark=otimes,
        style=ultra thick,
        ]
        coordinates {
        (50,7.407464779538410e-05)(100,1.747752187911350e-05)(200,2.934125327076634e-06)(400,4.046301042042377e-07)
        };
    
        \addplot[
        color=vermillion,
        mark=otimes,
        style=ultra thick,
        ]
        coordinates {
        (50,7.407464830230088e-05)(100,1.747752232155602e-05)(200,2.934125429187879e-06)(400,4.046301176938811e-07)
        };
    \legend{{} {$|\mathcal{T}_k^f| : |\mathcal{T}_k^s|$ = 1:1},{} {$|\mathcal{T}_k^f| : |\mathcal{T}_k^s|$ = 2:1},{} {$|\mathcal{T}_k^f| : |\mathcal{T}_k^s|$ = 4:1},{} {$|\mathcal{T}_k^f| : |\mathcal{T}_k^s|$ = 8:1},{} {$|\mathcal{T}_k^f| : |\mathcal{T}_k^s|$ = 16:1},}
    
    \addplot[
        color=gray,
        style=ultra thick,
        no markers,
        dashed
        ]
        coordinates {
            (50,0.00024407464)(400,4.7670828125000055e-07)
    };
    \node[] at (axis cs: 300, 2.5e-6) {\color{gray!60!black}\footnotesize $O(N^{-3})$};

    \end{axis}
    \end{tikzpicture}
    \caption{Convergence plot for the 2+1D heat and wave equation problem with solid source term with a finer fluid temporal mesh and uniform refinement in time}
    \label{fig:convergence_space_time_2d_solid_source_finer_fluid}
\end{figure}

{Using a finer fluid temporal mesh, we get the convergence plot in Figure \ref{fig:convergence_space_time_2d_solid_source_finer_fluid}. Therein, we observe that using a finer temporal mesh for the fluid in comparison to the solid does not improve the error in the quantity of interest. This was to be expected since we only have a non-zero right hand side in the solid domain and thus the solution behavior is mainly dictated by the wave equation.}

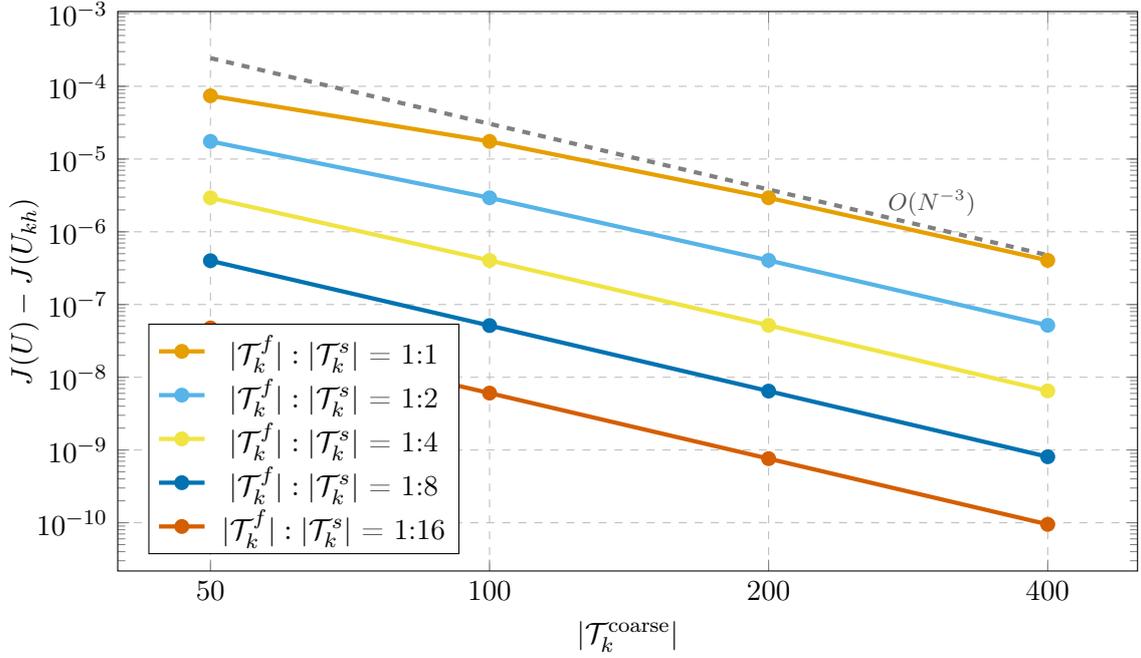
\begin{figure}[H]%
    \centering
    \begin{tikzpicture}
    \begin{axis}[
        xlabel={$|\mathcal{T}_k^{\text{coarse}}|$},
        ylabel={$J(U) - J(U_{kh})$},
        xmin=0.0, xmax=500,
        legend pos= south west,
        legend style={nodes={scale=1, transform shape}},
        xmajorgrids=true,
        ymajorgrids=true,
        grid style=dashed,
        xmode = log,
        ymode = log,
        width = 15cm, 
        height = 9cm,
        xtick = {50,100,200,400},
        xticklabels = {50,100,200,400},
    ]

    \addplot[
        color=orange,
        mark=otimes,
        style=ultra thick,
        ]
        coordinates {
        (50,7.407397938424647e-05)(100,1.747727249566005e-05)(200,2.934075263809842e-06)(400,4.046231868823792e-07)
        };
    
        \addplot[
        color=skyblue,
        mark=otimes,
        style=ultra thick,
        ]
        coordinates {
        (50,1.747497564090819e-05)(100,2.933719521883908e-06)(200,4.045765924470378e-07)(400,5.180853846311392e-08)
        };
    
        \addplot[
        color=yellow,
        mark=otimes,
        style=ultra thick,
        ]
        coordinates {
        (50,2.930395693736729e-06)(100,4.041892266948245e-07)(200,5.176211724489902e-08)(400,6.498447385349000e-09)
        };
    
        \addplot[
        color=blue,
        mark=otimes,
        style=ultra thick,
        ]
        coordinates {
        (50,4.005012657566686e-07)(100,5.136671143860308e-08)(200,6.452041335339338e-09)(400,8.069045410140363e-10)
        };
    
        \addplot[
        color=vermillion,
        mark=otimes,
        style=ultra thick,
        ]
        coordinates {
        (50,4.760579000223555e-08)(100,6.055398211503304e-09)(200,7.605001248970483e-10)(400,9.494325833338257e-11)
        };
    \legend{{} {$|\mathcal{T}_k^f| : |\mathcal{T}_k^s|$ = 1:1},{} {$|\mathcal{T}_k^f| : |\mathcal{T}_k^s|$ = 1:2},{} {$|\mathcal{T}_k^f| : |\mathcal{T}_k^s|$ = 1:4},{} {$|\mathcal{T}_k^f| : |\mathcal{T}_k^s|$ = 1:8},{} {$|\mathcal{T}_k^f| : |\mathcal{T}_k^s|$ = 1:16},}
    
    \addplot[
        color=gray,
        style=ultra thick,
        no markers,
        dashed
        ]
        coordinates {
            (50,0.00024407464)(400,4.7670828125000055e-07)
    };
    \node[] at (axis cs: 300, 2.5e-6) {\color{gray!60!black}\footnotesize $O(N^{-3})$};

    \end{axis}
    \end{tikzpicture}
    \caption{Convergence plot for the 2+1D heat and wave equation problem with solid source term with a finer solid temporal mesh and uniform refinement in time}
    \label{fig:convergence_space_time_2d_solid_source_finer_solid}
\end{figure}

{
Performing the same convergence test for a finer solid temporal mesh, the results are shown in Figure \ref{fig:convergence_space_time_2d_solid_source_finer_solid}. 
Here, we observe that using a finer temporal mesh for the solid in comparison to the fluid greatly reduces the error in the quantity of interest.
Instead of refining the temporal mesh uniformly for both fluid and solid, we can also refine the temporal mesh only for the solid and get similar results for our quantity of interest. 
}

\subsection{\auth{2+1D} Mandel's problem}
\label{sec_mandel_problem}
\begin{figure}[H]
    \begin{center}
    \begin{tikzpicture}[scale = 4, draw=black]
        \draw[draw=black, fill=black!5!white] (0,0) rectangle (3,0.6);
        \node (omega) at (1.5,0.3) {\Large{$\Omega$}};
        \node (left) at (-0.125,0.3) {$\Gamma_{\text{left}}$};
        \node (right) at (3.15,0.3) {$\Gamma_{\text{right}}$};
        \node (top) at (1.5,0.7) {$\Gamma_{\text{top}}$};
        \node (bottom) at (1.5,-0.1) {$\Gamma_{\text{bottom}}$};
    \end{tikzpicture}
    \caption{Domain for Mandel's problem}
    \label{fig:mandel_domain}
    \end{center}
\end{figure}
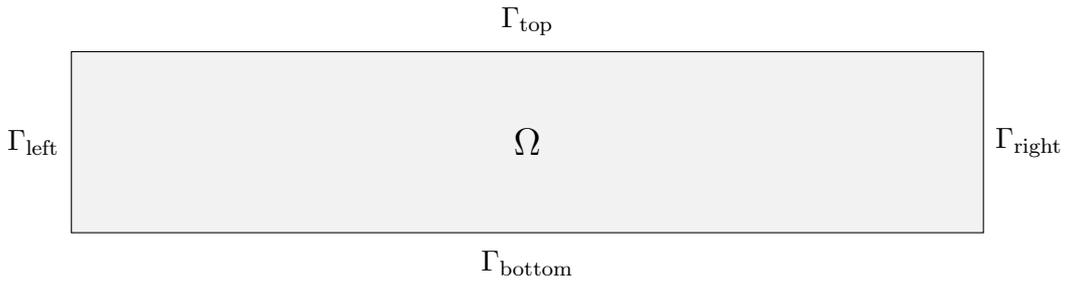

\auth{Next}, as an example for volume coupled problems, we consider Mandel's problem \cite{Mandel1953, Gai2004, Guzman2012, Girault2011, Liu04, Wick2020PFF, AbChCuDeRo96, Cheng88}, which is a benchmark from poroelasticity from Section \ref{sec_model_problem_poroelasticity},
and more recently also for nonlinear poroelasticity \cite{DuiMiWi22}.
This problem has also been solved by decoupled multirate schemes in \cite{Dean2006, Almani2016, Bause2017, Borregales2019}.
Therein, one observes the so-called Mandel-Cryer effect \cite{Cr63} of 
a non-monotonic pressure evolution: first increasing pressure, followed by decreasing pressure. 
Let $\Omega := (\SI{0}{\meter}, \SI{100}{\meter}) \times (\SI{0}{\meter}, \SI{20}{\meter})$ with boundaries as shown in Figure \ref{fig:mandel_domain}.
\noindent The initial and boundary conditions are given by
\begin{align*}
    p(\cdot, 0) &= p^0(\cdot) = 0 &&\text{in } \Omega \times \{0\}, \\
    u(\cdot, 0) &= u^0(\cdot) = 0 &&\text{in } \Omega \times \{0\}, \\
    \auth{\frac{K}{\nu}\nabla_x p\cdot n} &\auth{= 0} &&\auth{\text{on } \partial \Omega \setminus \Gamma_{\text{right}} \times I}, \quad\tag*{\auth{(No flow condition, homogeneous Neumann)}}\\
    \sigma(u) \cdot n &= -\bar{t}e_y &&\text{on } \Gamma_{\text{top}} \times I, \qquad\qquad\qquad\text{\auth{(Traction condition, inhomogeneous Neumann)}}\\
    p &= 0  &&\text{on } \Gamma_{\text{right}} \times I, \quad\tag*{\auth{(Constant zero pressure, homogeneous Dirichlet)}}\\
    \auth{\sigma(u) \cdot n} &\auth{= 0} &&\auth{\text{on } \Gamma_{\text{right}} \times I,}\quad\tag*{\auth{(Traction free, homogeneous Neumann)}}\\
    u_y &= 0 \quad\text{and}\quad \auth{\partial_y u_x = 0}  &&\text{on } \Gamma_{\text{bottom}} \times I,\qquad\qquad\text{\auth{(Confined conditions, mixed Dirichlet/Neumann)}} \\
    u_x &= 0  \quad\text{and}\quad \auth{\partial_x u_y = 0} &&\text{on } \Gamma_{\text{left}} \times I. \quad\tag*{\auth{(Confined conditions, mixed Dirichlet/Neumann)}}
\end{align*}
The parameters for Mandel's problem are summarized in Table \ref{tab:params_mandel}.
\begin{table}[H]
    \centering
    \begin{tabular}{ |p{5cm}||p{3cm}|}
         \hline
         Parameter & Value \\
         \hline
            M & \SI{1.75e7}{\pascal} \\
            c & 1/M \\
            $\alpha$ & \SI{1}{\pascal\metre} \\
            $\nu$ & \SI{1e-3}{\metre\squared\per\second} \\
            K & \SI{1e-13}{\metre\squared} \\
            $\rho$ & \SI{1}{\kilogram\per\metre\cubed} \\
            $\bar{t}$ & \SI{1e7}{\pascal\metre} \\
             $g$ & $0$ \\
            $q$ & $0$ \\ 
            $f$ & $0$ \\
            $\mu$ & \SI{1e8}{} \\
            $\lambda$ & $\frac{2}{3} \times 10^8$ \\
         \hline
    \end{tabular}
    \caption{Parameters in Mandel's problem} 
    \label{tab:params_mandel}
\end{table}
\noindent We employ a space-time finite element discretization and a Galerkin discretization with $\dG(0)$ in time for this numerical test. In space, we use quadratic finite elements for the displacement $u$ and linear finite elements for the pressure $p$. 

\begin{figure}[H]
     \centering
     \subfloat[pressure]{
    \includegraphics[width=14cm]{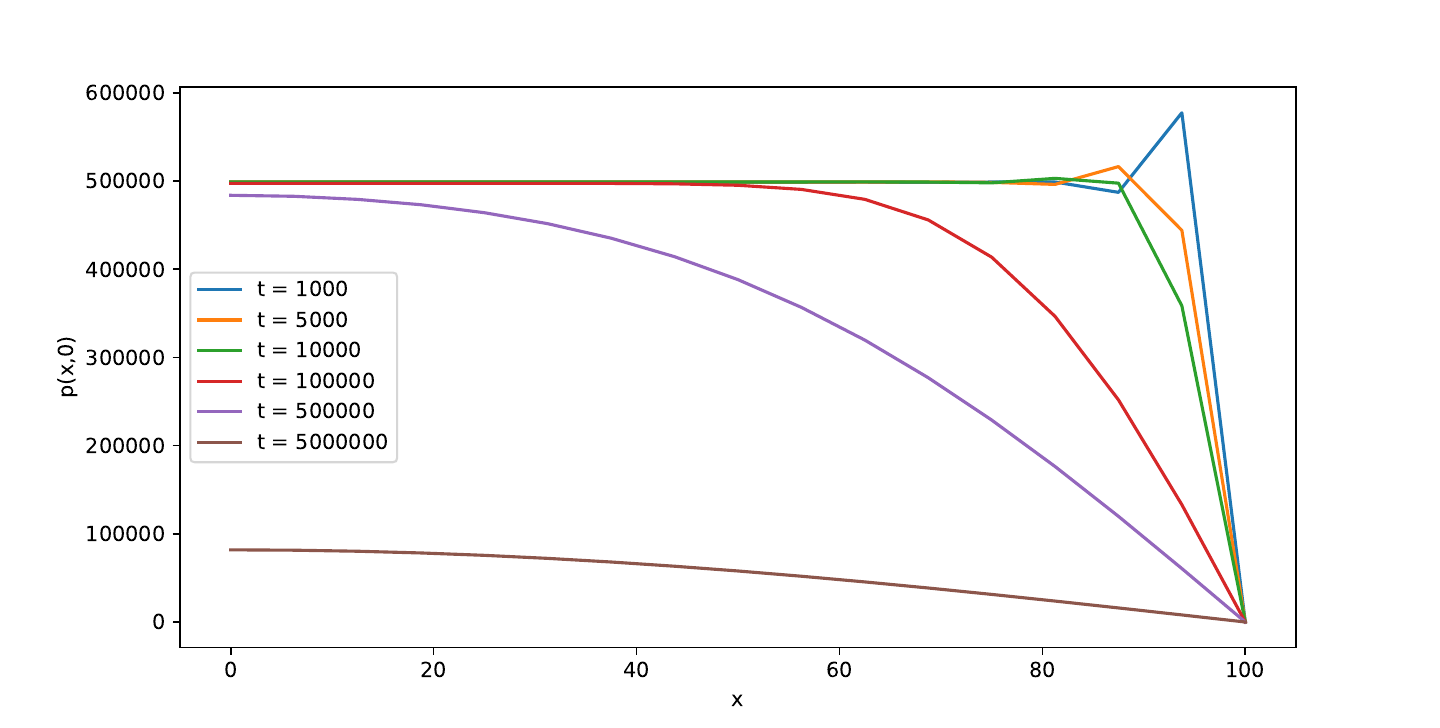}}%
    \\ 
     \subfloat[x-displacement]{
    \includegraphics[width=14cm]{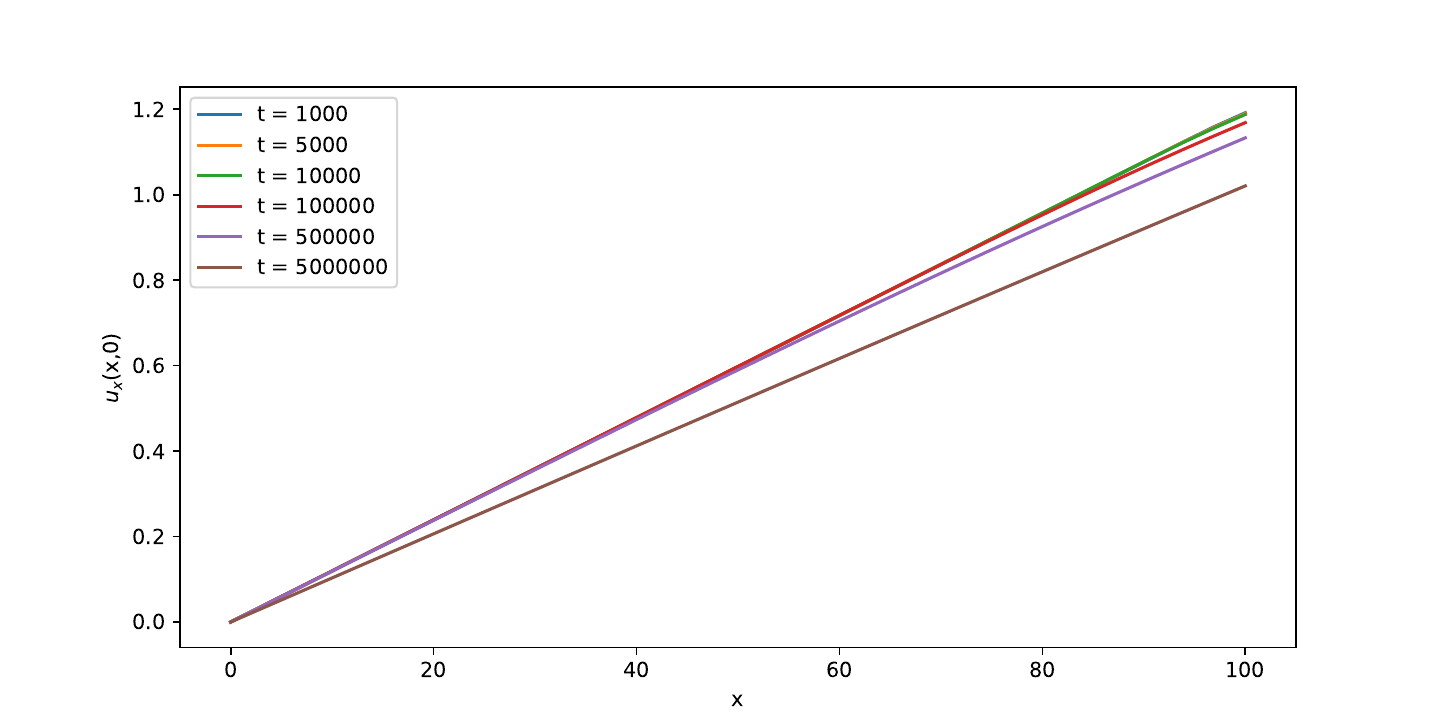}}%
    \caption{Reference solution on the bottom boundary ($\Gamma_{\text{bottom}}$) for the Mandel problem}
    \label{fig:reference_solution_Mandel}
\end{figure}

For this numerical test, we use a fixed spatial mesh with 16 spatial cells in the x-direction and 16 spatial cells in the y-direction.
As a reference solution, we solve the Mandel with $\dG(0)$ elements and $500,000$ temporal elements for displacement and pressure. The solution at the bottom boundary is shown in Figure \ref{fig:reference_solution_Mandel}. 

To validate the efficiency of our proposed multirate methodology, we use the quantity of interest
\begin{align*}
    J(U_{kh}) &= \int_I \int_{\Gamma_{\text{bottom}}} p_{kh} (x,t)\ \mathrm{d}(x,t),
\end{align*}
and compare with the reference value $J(U) := 8.718831\cdot10^{13}$ from the simulation with $500,000$ temporal elements.

For the convergence studies, we employ an initial coarse temporal mesh with $|\mathcal{T}_k^{\text{coarse}}| = 1250$ temporal elements and then uniformly refine the displacement or pressure temporal mesh up to 4 times. This means that e.g. for each displacement temporal element we have 16 pressure temporal elements.

\begin{figure}[H]%
    \centering
    \begin{tikzpicture}
    \begin{axis}[
        xlabel={$|\mathcal{T}_k^{\text{coarse}}|$},
        ylabel={$J(U) - J(U_{kh})$},
        xmin=0.0, xmax=13000,
        legend pos= north east, 
        legend style={nodes={scale=1, transform shape}},
        xmajorgrids=true,
        ymajorgrids=true,
        grid style=dashed,
        xmode = log,
        ymode = log,
        width = 15cm, 
        height = 9cm,
        xtick = {1250,2500,5000,10000},
        xticklabels = {1250,2500,5000,10000},
        ytick = {2500000000, 5000000000,10000000000, 20000000000},
        yticklabels = {$2.5 \cdot 10^9$, $5 \cdot 10^9$, $10^{10}$, $2 \cdot 10^{10}$},
    ]
    
    \addplot[
    color=orange,
	mark=otimes,
    style=ultra thick,
    ]
    coordinates {
	(1250,2.125187231104688e+10)(2500,1.060065598395312e+10)(5000,5.274034088875000e+09)(10000,2.610469713203125e+09)
	};

    \addplot[
    color=skyblue,
	mark=otimes,
    style=ultra thick,
    ]
    coordinates {
	(1250,2.125187231625000e+10)(2500,1.060065597284375e+10)(5000,5.274034075796875e+09)(10000,2.610469766437500e+09)
	};

    \addplot[
    color=yellow,
	mark=otimes,
    style=ultra thick,
    ]
    coordinates {
	(1250,2.125187231475000e+10)(2500,1.060065598182812e+10)(5000,5.274034079031250e+09)(10000,2.610469722015625e+09)
	};

    \addplot[
    color=blue,
	mark=otimes,
    style=ultra thick,
    ]
    coordinates {
	(1250,2.125187231915625e+10)(2500,1.060065604120312e+10)(5000,5.274034093562500e+09)(10000,2.610469805406250e+09)
	};

    \addplot[
    color=vermillion,
	mark=otimes,
    style=ultra thick,
    ]
    coordinates {
	(1250,2.125187234764062e+10)(2500,1.060065597140625e+10)(5000,5.274034161640625e+09)(10000,2.610470049328125e+09)
	};

     \addplot[
    color=gray,
    style=ultra thick,
    no markers,
    dashed
    ]
    coordinates {
    (1250,2.5e+10)(10000,3.125e+9) 
    };
    \node[] at (axis cs: 3000, 1.25e10) {\color{gray!60!black}\footnotesize $O(N^{-1})$};

    \legend{{} {$|\mathcal{T}_k^u| : |\mathcal{T}_k^p|$ = 1:1},{} {$|\mathcal{T}_k^u| : |\mathcal{T}_k^p|$ = 2:1},{} {$|\mathcal{T}_k^u| : |\mathcal{T}_k^p|$ = 4:1},{} {$|\mathcal{T}_k^u| : |\mathcal{T}_k^p|$ = 8:1},{} {$|\mathcal{T}_k^u| : |\mathcal{T}_k^p|$ = 16:1},}
    \end{axis}
    \end{tikzpicture}
    \caption{Convergence plot for the Mandel problem with a finer displacement temporal mesh and uniform refinement in time}
    \label{fig:convergence_space_time_Mandel_finer_displacement}
\end{figure}
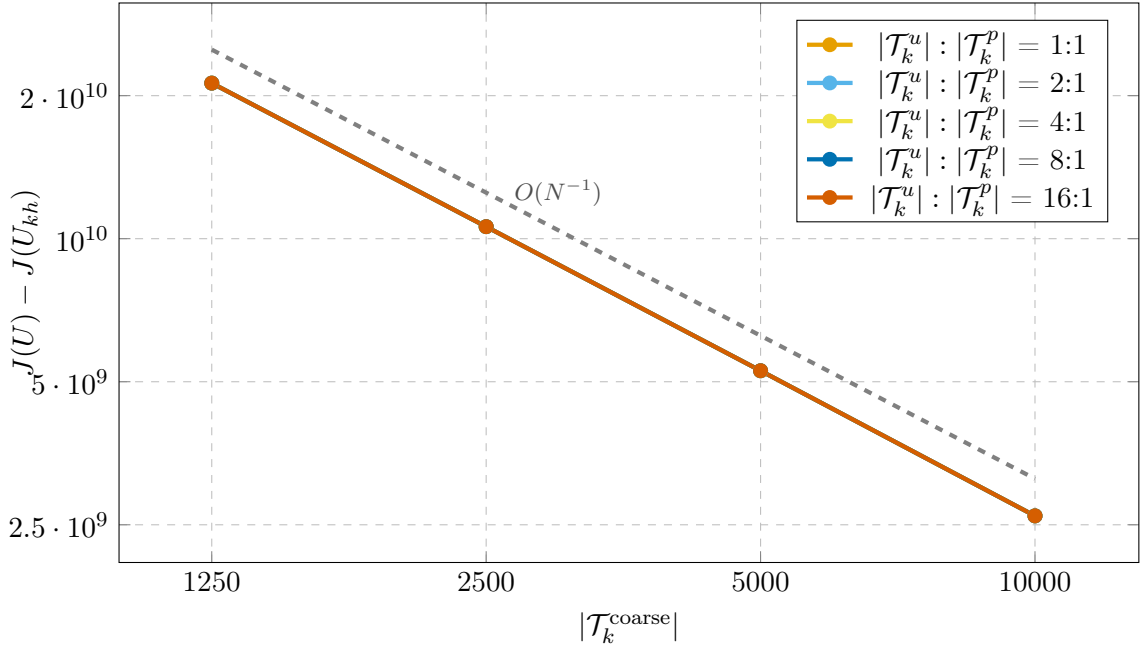

Using a finer displacement temporal mesh, we get the convergence plot in Figure \ref{fig:convergence_space_time_Mandel_finer_displacement}. 
Here, we see that using a finer temporal mesh for the displacement in comparison to the pressure does not improve the error in the quantity of interest and the error is identical for all temporal meshes in the first five significant digits. This is not surprising, since the constitutive equation for the displacement is a quasi-stationary problem and the quantity of interest only measures a pressure boundary integral over time. Therefore, one can expect only refinement of the temporal mesh of the pressure to lead to an error reduction.

\begin{figure}[H]%
    \centering
    \begin{tikzpicture}
    \begin{axis}[
        xlabel={$|\mathcal{T}_k^{\text{coarse}}|$},
        ylabel={$J(U) - J(U_{kh})$},
        xmin=0.0, xmax=13000,
        legend pos= south west,
        legend style={nodes={scale=1, transform shape}},
        xmajorgrids=true,
        ymajorgrids=true,
        grid style=dashed,
        xmode = log,
        ymode = log,
        width = 15cm, 
        height = 9cm,
        xtick = {1250,2500,5000,10000},
        xticklabels = {1250,2500,5000,10000},
        ytick = {100000000, 1000000000,10000000000},
        yticklabels = {$10^8$, $10^9$, $10^{10}$},
    ]

        \addplot[
    color=orange,
	mark=otimes,
    style=ultra thick,
    ]
    coordinates {
	(1250,2.125187231104688e+10)(2500,1.060065598395312e+10)(5000,5.274034088875000e+09)(10000,2.610469713203125e+09)
	};

    \addplot[
    color=skyblue,
	mark=otimes,
    style=ultra thick,
    ]
    coordinates {
	(1250,1.060089025304688e+10)(2500,5.274092728156250e+09)(5000,2.610484525750000e+09)(10000,1.278627673093750e+09)
	};

    \addplot[
    color=yellow,
	mark=otimes,
    style=ultra thick,
    ]
    coordinates {
	(1250,5.274327089890625e+09)(2500,2.610543177890625e+09)(5000,1.278642449968750e+09)(10000,6.126900334062500e+08)
	};

    \addplot[
    color=blue,
	mark=otimes,
    style=ultra thick,
    ]
    coordinates {
	(1250,2.610777580859375e+09)(2500,1.278701094265625e+09)(5000,6.127047630000000e+08)(10000,2.797169916093750e+08)
	};

    \addplot[
    color=vermillion,
	mark=otimes,
    style=ultra thick,
    ]
    coordinates {
	(1250,1.278935544843750e+09)(2500,6.127634227968750e+08)(5000,2.797317876718750e+08)(10000,1.132293052812500e+08)
	};
    \legend{{} {$|\mathcal{T}_k^u| : |\mathcal{T}_k^p|$ = 1:1},{} {$|\mathcal{T}_k^u| : |\mathcal{T}_k^p|$ = 1:2},{} {$|\mathcal{T}_k^u| : |\mathcal{T}_k^p|$ = 1:4},{} {$|\mathcal{T}_k^u| : |\mathcal{T}_k^p|$ = 1:8},{} {$|\mathcal{T}_k^u| : |\mathcal{T}_k^p|$ = 1:16},}

    \addplot[
    color=gray,
    style=ultra thick,
    no markers,
    dashed
    ]
    coordinates {
    (1250,2.5e+10)(10000,3.125e+9) 
    };
    \node[] at (axis cs: 3000, 1.45e10) {\color{gray!60!black}\footnotesize $O(N^{-1})$};
    
    \end{axis}
    \end{tikzpicture}
    \caption{Convergence plot for the Mandel problem with a finer pressure temporal mesh and uniform refinement in time}
    \label{fig:convergence_space_time_mandel_finer_pressure}
\end{figure}
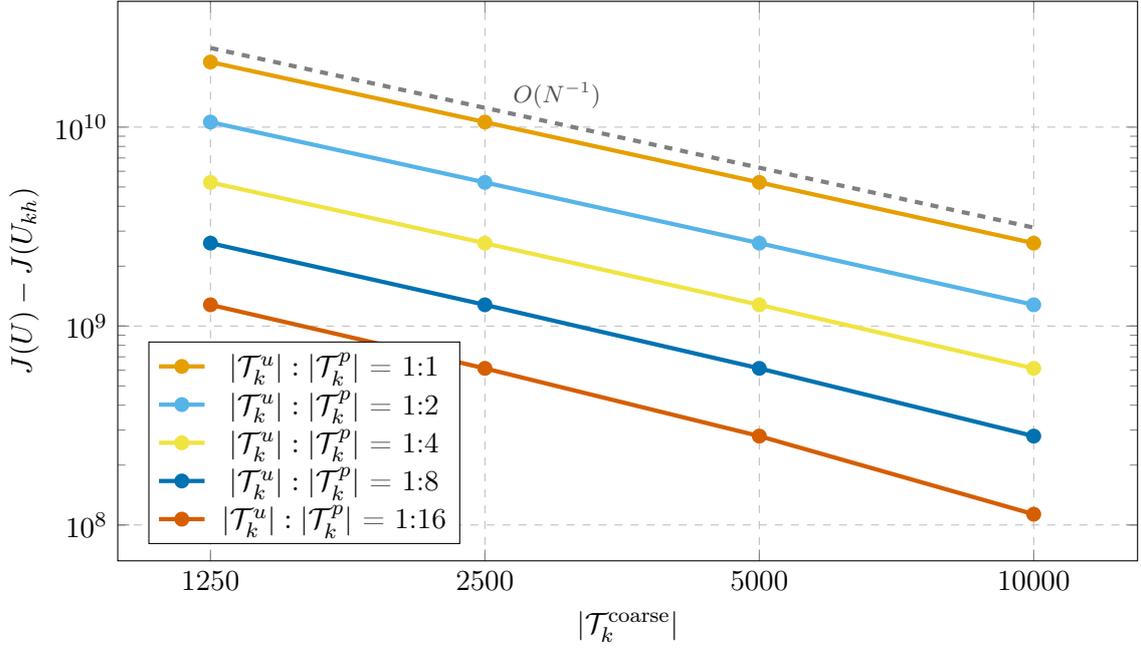

In Figure \ref{fig:convergence_space_time_mandel_finer_pressure}, we run the same convergence test for a finer temporal mesh for the pressure.
Now, it pays off using a finer temporal mesh for the pressure than for the displacement. Not only do we have linear convergence of the error under uniform refinement of both temporal meshes, but also for the refinement ratio of pressure to displacement. 
More concretely, the error roughly halves when going from 1:1 to 1:2, going from 1:2 to 1:4, etc. 
This shows that instead of using conventional time stepping schemes, which rely on the same temporal mesh for displacement and pressure, we can get the same error reduction by only increasing the number of space-time degrees of freedom for the pressure and using a fixed number of space-time degrees of freedom for the displacement.
This multirate scheme is particularly cost-effective when the displacement requires high-order spatial finite elements than the pressure, e.g. as in the case of Taylor-Hood elements. 

\subsection{3+1D footing problem}
\label{sec_footing_problem}
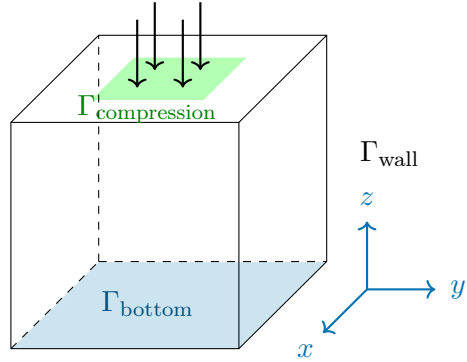
\begin{figure}[H]
    \begin{center}
    \begin{tikzpicture}[scale=3, draw=black]
        \draw[draw=white,fill=blue!20!white] (0,0,0) -- (1,0,0) -- (1,0,-1) -- (0,0,-1) -- (0,0,0);
        \draw[draw=white,fill=green!30!white] (0.25,1,-0.25) -- (0.25,1,-0.75) -- (0.75,1,-0.75) -- (0.75,1,-0.25) -- (0.25,1,-0.25);
        \draw[draw=black] (0,0,0) -- (1,0,0) -- (1,1,0) -- (0,1,0) -- (0,0,0);
        \draw[draw=black] (1,0,0) -- (1,0,-1) -- (1,1,-1) -- (0,1,-1) -- (0,1,0);
        \draw[draw=black] (1,1,0) -- (1,1,-1);
        \draw[draw=black, dashed] (0,0,0) -- (0,0,-1);
        \draw[draw=black, dashed] (0,0,-1) -- (1,0,-1);
        \draw[draw=black, dashed] (0,0,-1) -- (0,1,-1);
        \node (bottom) at (0.6,0.2,0) {\color{blue!70!black}$\Gamma_{\text{bottom}}$};
        \node (compression) at (0.29,0.75,-0.8) {\color{green!50!black}$\Gamma_{\text{compression}}$};
        \node (wall) at (1.2,0.4,-1.2) {$\Gamma_{\text{wall}}$};
        \draw[draw=black, ->, thick]
        (0.4, 1.3, -0.4) -- (0.4, 1.0, -0.4);
        \draw[draw=black, ->, thick]
        (0.4, 1.3, -0.6) -- (0.4, 1.0, -0.6);
        \draw[draw=black, ->, thick]
        (0.6, 1.3, -0.4) -- (0.6, 1.0, -0.4);
        \draw[draw=black, ->, thick]
        (0.6, 1.3, -0.6) -- (0.6, 1.0, -0.6);
        \draw[draw=blue, ->, thick]
        (1.1,-0.2,-1.2) -- (1.4,-0.2,-1.2);
        \node(x) at (1.5, -0.2, -1.2) {\color{blue}$y$};
        \draw[draw=blue, ->, thick]
        (1.1,-0.2,-1.2) -- (1.1,0.1,-1.2);
        \node(z) at (1.1, 0.2, -1.2) {\color{blue}$z$};
        \draw[draw=blue, ->, thick]
        (1.1,-0.2,-1.2) -- (1.1,-0.2,-0.7);
        \node(y) at (1.1, -0.2, -0.5) {\color{blue}$x$};
    \end{tikzpicture}
    \caption{Domain for 3D footing problem.}
    \label{fig:footing_domain}
    \end{center}
\end{figure}

{
In this final numerical example, a three-dimensional footing problem inspired by \cite{Gaspar2008} is studied as another test for poroelasticity, i.e. for volume coupling. Let $\Omega := (\SI{-32}{\meter}, \SI{32}{\meter}) \times (\SI{-32}{\meter}, \SI{32}{\meter}) \times (\SI{0}{\meter}, \SI{64}{\meter})$ and $I := (\SI{0}{\second}, \SI{5000000}{\second})$ with boundaries as shown in Figure~\ref{fig:footing_domain}.
}
{
The initial and boundary conditions are given by
\begin{align*}
    p(0) &= p^0 = 0 &&\text{in } \Omega \times \{0\}, \\
    u(0) &= u^0 = 0 &&\text{in } \Omega \times \{0\}, \\
     \frac{K}{\nu}\nabla_x p\cdot n &= 0 &&\text{on } \partial \Omega \setminus \Gamma_{\text{bottom}} \times I, \quad\tag*{(No flow condition, homogeneous Neumann)}\\
    \sigma(u) \cdot n &= -\bar{t}e_z &&\text{on } \Gamma_{\text{compression}} \times I,\quad\tag*{(Traction condition, inhomogeneous Neumann)} \\
    \sigma(u) \cdot n &= 0 &&\text{on } \Gamma_{\text{top}} \setminus \Gamma_{\text{compression}}\times I,\quad\tag*{(Traction-free, homogeneous Neumann)}\\ 
    p &= 0  &&\text{on } \Gamma_{\text{bottom}} \times I,\quad\tag*{(Constant zero pressure, homogeneous Dirichlet)} \\
    u &= 0  &&\text{on } \Gamma_{\text{bottom}} \times I,\quad\tag*{(Fixed displacements, homogeneous Dirichlet)} \\
    \sigma(u) \cdot n &= 0 &&\text{on } \Gamma_{\text{wall}} \times I, \quad\tag*{(Traction-free, homogeneous Neumann)}
\end{align*}
The material parameters as listed in Table~\ref{tab:params_mandel} 
(Section \ref{sec_mandel_problem}) are used again.
}

\indent {
We employ a space-time finite element discretization and a Galerkin discretization with $\dG(0)$ in time for this numerical test. In space, we use Taylor-Hood elements, namely quadratic finite elements in space for the displacement $u_{kh}$ and linear finite elements for the pressure $p_{kh}$ are employed. 
The spatial mesh is fixed, it comprises 8 spatial cells in each direction i.e., an isotropic mesh with $14,739$ DoFs for displacement and $729$ DoFs for pressure.
As a reference solution, we solve the 3D footing problem with $\dG(0)$ in time and $50,000$ temporal elements for displacement and pressure. As a quantity of interest, we use the time-integrated pressure acting at the compression boundary i.e.,
\begin{align*}
    J(U_{kh}) := \int_I \int_{\Gamma_{\text{compression}}} p_{kh}\ \mathrm{d}x\ \mathrm{d}t,
\end{align*}
and compare with the reference value $J(U) := 1.207538 \cdot 10^{14}$ from the simulation with $50,000$ temporal elements.
}

{
For the convergence studies, we employ an initial coarse temporal mesh with $|\mathcal{T}_k^{\text{coarse}}| = 125$ temporal elements and then uniformly refine the pressure temporal mesh up to 4 times. This means that e.g. for each displacement temporal element we have 16 pressure temporal elements. For this numerical example, we only consider the case of a finer pressure temporal mesh than displacement temporal mesh, since this was the case where we observed the largest error reduction for Mandel's benchmark problem.
Moreover, for more displacement temporal elements efficient iterative solvers should be used, whereas we use a direct solver for the linear system of equations in this numerical example.
}

\begin{figure}[H]%
    \centering
    \begin{tikzpicture}
    \begin{axis}[
        xlabel={$|\mathcal{T}_k^{\text{coarse}}|$},
        ylabel={$J(U) - J(U_{kh})$},
        xmin=0.0, xmax=1300,
        legend pos= south west,
        legend style={nodes={scale=1, transform shape}},
        xmajorgrids=true,
        ymajorgrids=true,
        grid style=dashed,
        xmode = log,
        ymode = log,
        width = 15cm, 
        height = 9cm,
        xtick = {125,250,500,1000},
        xticklabels = {125,250,500,1000},
        ytick = {1000000000,10000000000,100000000000},
        yticklabels = {$10^9$, $10^{10}$, $10^{11}$},
    ]

    \addplot[
        color=orange,
        mark=otimes,
        style=ultra thick,
        ]
        coordinates {
        (125,1.819371106172812e+11)(250,9.059873332823438e+10)(500,4.503571171067188e+10)(1000,2.228139565484375e+10)
        };
    
        \addplot[
        color=skyblue,
        mark=otimes,
        style=ultra thick,
        ]
        coordinates {
        (125,9.066991197353125e+10)(250,4.505305830596875e+10)(500,2.228564218993750e+10)(1000,1.091216507853125e+10)
        };
    
        \addplot[
        color=yellow,
        mark=otimes,
        style=ultra thick,
        ]
        coordinates {
        (125,4.512313284654688e+10)(250,2.230271810007812e+10)(500,1.091635903406250e+10)(1000,5.229007964609375e+09)
        };
    
        \addplot[
        color=blue,
        mark=otimes,
        style=ultra thick,
        ]
        coordinates {
        (125,2.237203293367188e+10)(250,1.093325900189062e+10)(500,5.233169960234375e+09)(1000,2.387796507906250e+09)
        };
    
        \addplot[
        color=vermillion,
        mark=otimes,
        style=ultra thick,
        ]
        coordinates {
        (125,1.100212429485938e+10)(250,5.249967747953125e+09)(500,2.391940361812500e+09)(1000,9.672829138906250e+08)
        };
    \legend{{} {$|\mathcal{T}_k^u| : |\mathcal{T}_k^p|$ = 1:1},{} {$|\mathcal{T}_k^u| : |\mathcal{T}_k^p|$ = 1:2},{} {$|\mathcal{T}_k^u| : |\mathcal{T}_k^p|$ = 1:4},{} {$|\mathcal{T}_k^u| : |\mathcal{T}_k^p|$ = 1:8},{} {$|\mathcal{T}_k^u| : |\mathcal{T}_k^p|$ = 1:16},}
    
    \addplot[
    color=gray,
    style=ultra thick,
    no markers,
    dashed
    ]
    coordinates {
    (125,211940000000.0)(1000,26492499999.999985) 
    };
    \node[] at (axis cs: 300, 1.25e11) {\color{gray!60!black}\footnotesize $O(N^{-1})$};
    
    \end{axis}
    \end{tikzpicture}
    \caption{Convergence plot for the 3+1D footing problem with a finer pressure temporal mesh and uniform refinement in time}
    \label{fig:convergence_space_time_footing_finer_pressure}
\end{figure}
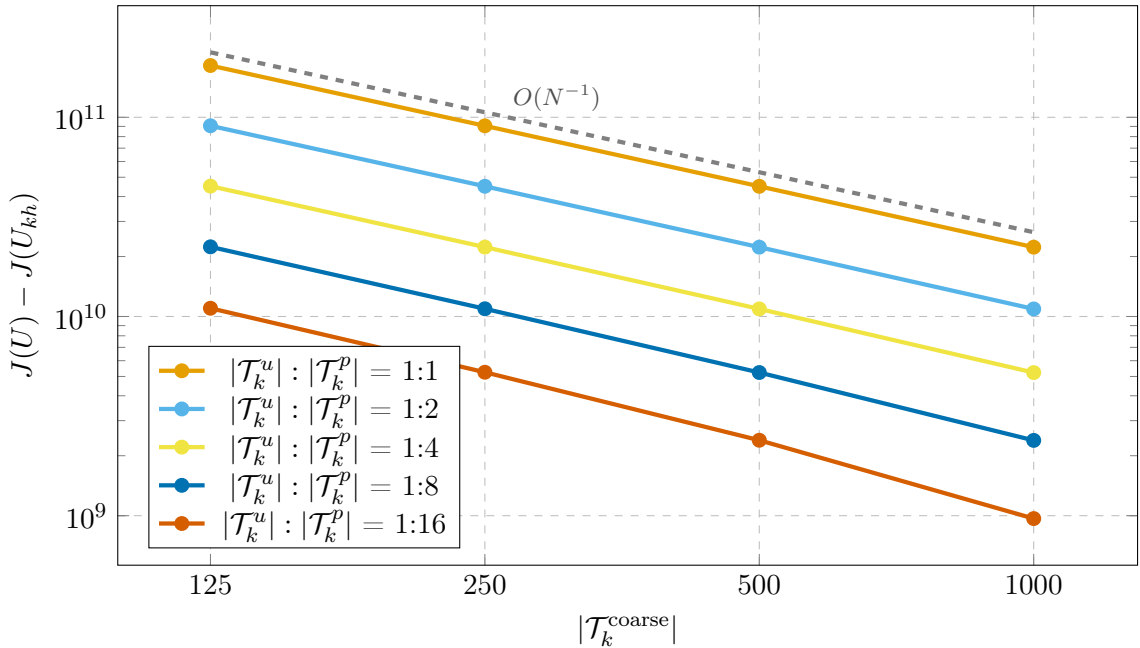

{
Using a finer pressure temporal mesh, we get the convergence plot in Figure \ref{fig:convergence_space_time_footing_finer_pressure}.
We make similar observations as for Mandel's problem, which, on the one hand, 
was to be expected since we use the same equations. 
On the other hand, the extension from 2D (in space) to 3D is a major step as all formulations 
need to be extended correspondingly, which requires all the programming code (including 
implementations and debugging) to be modified to 3D. Thus, our findings indicate the correctness 
of our implementations and robustness of our approach independently of the spatial dimension.
Again in more detail, when refining both temporal meshes uniformly, we observe linear convergence of the error in the quantity of interest.
Similarly, when just refining the pressure temporal mesh, we observe linear convergence of the error in the quantity of interest.
This motivates the usage of monolithic temporal multirate schemes especially for poroelasticity, since we can get the same error reduction by only increasing the number of temporal elements for pressure and using a fixed temporal mesh for displacement in our numerical examples.
}

\section{Conclusions and outlook}
\label{sec_conclusions}
In this work, we proposed a monolithic space-time multirate framework for 
interface-coupled and volume-coupled problems.
One first objective was a mathematical abstract formalism that allows to consider both interface-coupled and volume-coupled problems.
As such model problems, 
an interface-coupled heat-wave system and a volume-coupled poro-elasticity problem 
were taken. The function spaces and tensor-product Galerkin finite element 
discretization were designed to treat the multirate idea in a monolithic fashion. 
Here, classical continuous finite elements were employed for the spatial 
discretization, while the temporal approximation was based on discontinuous finite elements. The concepts and numerical realization were discussed in great detail in Section \ref{sec_methodology}.
{Five} numerical tests {for different spatial dimensions} 
were conducted in which the multirate schemes were computationally
analyzed in detail. Here, goal-functionals were adopted in order to study specific 
quantities of interest. Overall, the performances were excellent.
So far, the different time meshes were given a priori. A future extension is 
to determine the time meshes by an error estimator; see also \cite{Logg2003a} for early work in this direction
applied to ODEs. We have already used 
goal functionals as quantities of interest in the current paper. Consequently,
goal-oriented a posteriori error estimates will be a natural choice to
obtain such time meshes. 
Finally, our abstract formalism, corresponding algorithms, and implementation were designed in such a way that they allow 
for future extensions to nonlinear problems.

\appendix

\section{Monolithic multirate backward Euler scheme}
\label{sec_monolithic_multirate_backward_euler}

In this section, we will discuss how time marching approaches can be used and extended to solve multirate in time formulations of multiphysics problems.
We will explain the methodology with an example of the backward Euler scheme: 
Find $U_{kh} \in X_k^{\dG(0)}\left((\mathcal{T}_k^1, \mathcal{T}_k^2), V_h(\mathcal{T}_h)\right)$ such that
\begin{align*}
    \tilde{A}(U_{kh})(\Phi_{kh})  = \tilde{F}(\Phi_{kh}) \qquad \forall \Phi_{kh} \in X_k^{\dG(0)}\left((\mathcal{T}_k^1, \mathcal{T}_k^2), V_h(\mathcal{T}_h)\right).
\end{align*}
From the definition of our time mesh structure, we have
$
\mathcal{T}_k^{\text{coarse}} \subset \mathcal{T}_k^1, \mathcal{T}_k^2 \subset \mathcal{T}_k^{\text{fine}}
$. We use this hierarchy to establish a method that borrows from both standard time marching methodology as well as space-time approach. Let us assume some fixed $I_{n-1}, I_{n}, I_{n + 1} \in \mathcal{T}_k^{\text{coarse}}$.
Since the temporal meshes $\mathcal{T}_k^1$ and $\mathcal{T}_k^2$ originate from the common coarse mesh through adaptive refinement, for each coarse element $I_n$ we can find elements $I_{n_1}^1, ..., I_{n_{N_n^1}}^1 \in \mathcal{T}_k^1$ and $I_{n_1}^2, ..., I_{n_{N_n^2}}^2 \in \mathcal{T}_k^2$ such that
\begin{equation}
\bar{I}_n = \bigcup_{N = 1}^{N_n^1} \bar{I}_{n_N}^1 = \bigcup_{N = 1}^{N_n^2} \bar{I}_{n_N}^2.
\label{time_partition}
\end{equation}
Our approach is based on simultaneously solving all of the equations given by the partition~(\ref{time_partition}) while using the last available solution from $I_{n-1}$ with respect to either $\mathcal{T}_k^1$ or $\mathcal{T}_k^2$ as the initial value. The time interpolation necessary to compute coupling conditions at temporal hanging nodes, i.e.  degrees of freedom that are contained in $\mathcal{T}_k^1$ and not in $\mathcal{T}_k^2$ or vice versa, is given by 
\begin{equation*}
\bar{u}_i \Big{|}_{I_n} \coloneqq \frac{1}{|I_n|}\int_{I_n} u_i \ \mathrm{d}t\hspace*{0.5 cm} \text{for }i\in\{1,2\}.
\end{equation*}

\section{Derivation of a time-stepping scheme for the Mandel problem}
\label{sec_timestepping_mandel}

In the following, we will derive by hand the time-stepping scheme for the Mandel problem with one displacement temporal element and two pressure temporal elements. Other multirate time-stepping schemes can be derived in a similar fashion.
For this we take the variational formulation (\ref{eq:variational_form_mandel}) of the Mandel problem, remove the terms in the right hand side function which are zero for the Mandel problem and only consider a time step $I_m = (0, k)$. For the displacement, we have only the temporal element $(0,k)$ and for the pressure we have the temporal elements $(0, \frac{k}{2})$ and $(\frac{k}{2},k)$.
For $\dG(0)$ in time we then have the ansatz
\begin{align*}
    u(t) &= u_1 \cdot \underbrace{\chi_{(0,k)}(t)}_{=: \phi_k^{u, (1)}}, \\
    p(t) &= p_1 \cdot \underbrace{\chi_{(0,\frac{k}{2})}(t)}_{=: \phi_k^{p, (1)}} + p_2 \cdot \underbrace{\chi_{(\frac{k}{2},k)}(t)}_{=: \phi_k^{p, (2)}},
\end{align*}
where $\chi_{I_k}(\cdot)$ denotes the temporal indicator function of the temporal element $I_k$.
We introduce the spatial matrices (and vector)
\begin{align*}
    K^p &= \left\{ \frac{K}{\nu}(\nabla_x \phi_h^{p,(j)},\nabla_x \phi_h^{p,(i)}) \right\}_{i,j = 1}^{\# \text{DoFs}(\mathcal{T}_h^p)}, \\
    M^p &= \left\{ c(\phi_h^{p,(j)}, \phi_h^{p,(i)}) \right\}_{i,j = 1}^{\# \text{DoFs}(\mathcal{T}_h^p)}, \\
    \Sigma^u &= \left\{ (\sigma(\phi_h^{u,(j)}),\nabla_x\phi_h^{u,(i)}) \right\}_{i,j = 1}^{\# \text{DoFs}(\mathcal{T}_h^u)}, \\
    F^u &= \left\{ \langle -\bar{t}e_y, \phi_h^{u, (i)}\rangle_{ \Gamma_{\text{top}}} \right\}_{i = 1}^{\# \text{DoFs}(\mathcal{T}_h^u)}, \\
    B^{up} &= \left\{ -\alpha(\phi_h^{p,(j)}I,\nabla_x \phi_h^{u,(i)}) + \alpha\langle \phi_h^{p,(j)}n, \phi_h^{u,(i)}\rangle_{ \Gamma_{\text{top}}} \right\}_{i,j = 1}^{\# \text{DoFs}(\mathcal{T}_h^u), \# \text{DoFs}(\mathcal{T}_h^p)}, \\
    B^{pu} &= \left\{ \alpha(\nabla_x \cdot \phi_h^{u,(j)},\phi_h^{p,(i)})\right\}_{i,j = 1}^{\# \text{DoFs}(\mathcal{T}_h^p), \# \text{DoFs}(\mathcal{T}_h^u)}.
\end{align*}
Using these spatial matrices, the time stepping scheme for one displacement temporal element and two pressure temporal elements is given by
\begin{align*}
    &(\phi_k^{u, (1)}, \phi_k^{u, (1)})_{I_m} \Sigma^u u_1 + \left\{(\phi_k^{p, (j)}, \phi_k^{u, (1)})_{I_m}\right\}_{j=1}^2 \otimes B^{up} \begin{pmatrix}
        p_1 \\ p_2
    \end{pmatrix} =  (1, \phi_k^{u, (1)})_{I_m} F^u, \\
    &\left\{(\phi_k^{p, (j)}, \phi_k^{p, (i)})_{I_m}\right\}_{i,j=1}^2 \otimes K^{p} \begin{pmatrix}
        p_1 \\ p_2
    \end{pmatrix}  + 
    \begin{pmatrix}
        1 & 0 \\
        -1 & 1
    \end{pmatrix}
    \otimes M^p \begin{pmatrix}
        p_1 \\ p_2
    \end{pmatrix} \\
    &+ \begin{pmatrix}
        1 & 0
    \end{pmatrix}
    \otimes B^{pu} u_1
    = \begin{pmatrix}
        1 & 0
    \end{pmatrix}
    \otimes B^{pu} u_0
    + \begin{pmatrix}
        1 & 0
    \end{pmatrix}
    \otimes M^p p_0.
\end{align*}
Now we evaluate all temporal integrals from above
\begin{align*}
    (\phi_k^{u, (1)}, \phi_k^{u, (1)})_{I_m} &= k, \\
    (\phi_k^{p, (i)}, \phi_k^{p, (j)})_{I_m} &= \frac{k}{2}\delta_{ij}, \\
    (\phi_k^{p, (j)}, \phi_k^{u, (1)})_{I_m} &= \frac{k}{2},
\end{align*}
and arrive at the block linear system
\begin{align*}
    \begin{pmatrix}
        k\Sigma^u & \frac{k}{2}B^{up} & \frac{k}{2}B^{up} \\
        B^{pu} & \frac{k}{2}K^p + M^p & 0 \\
        0 & -M^p & \frac{k}{2}K^p + M^p
    \end{pmatrix}
    \begin{pmatrix}
        u_1 \\
        p_1 \\
        p_2
    \end{pmatrix}
    = 
    \begin{pmatrix}
        kF^u \\
        B^{pu}u_0 + M^p p_0 \\
        0
    \end{pmatrix}.
\end{align*}
In the space-time methodology from Section \ref{sec_temporal_multirate_tp_st_fem} the matrix 
has the same content, but the coupling blocks 
\begin{align*}
    \begin{pmatrix}
        B^{pu} \\ 0
    \end{pmatrix} \qquad \text{and} \qquad 
    \begin{pmatrix}
        \frac{k}{2}B^{up} & \frac{k}{2}B^{up}
    \end{pmatrix}
\end{align*}
are being computed by first assembling the temporal integrals with the finer temporal basis, i.e. we have
\begin{align*}
    \begin{pmatrix}
        B^{pu} & 0 \\ -B^{pu} & B^{pu}
    \end{pmatrix} \qquad \text{and} \qquad 
    \begin{pmatrix}
        \frac{k}{2}B^{up} & 0 \\ 0 & \frac{k}{2}B^{up}
    \end{pmatrix}.
\end{align*}
Using the restriction matrix
\begin{align*}
    R = \begin{pmatrix}
        1 & 1
    \end{pmatrix},
\end{align*}
we then have
\begin{align*}
    \begin{pmatrix}
        B^{pu} & 0 \\ -B^{pu} & B^{pu}
    \end{pmatrix} \begin{pmatrix}
        1 \\ 1
    \end{pmatrix} = \begin{pmatrix}
        B^{pu} \\ 0
    \end{pmatrix} \qquad \text{and} \qquad 
    \begin{pmatrix}
        1 & 1
    \end{pmatrix}
    \begin{pmatrix}
        \frac{k}{2}B^{up} & 0 \\ 0 & \frac{k}{2}B^{up}
    \end{pmatrix} = \begin{pmatrix}
        \frac{k}{2}B^{up} & \frac{k}{2}B^{up}
    \end{pmatrix}.
\end{align*}

\section*{Acknowledgements}
JR acknowledges the funding of the German Research Foundation (DFG) within the framework of the International Research Training Group on  Computational Mechanics Techniques in High Dimensions GRK 2657 under Grant Number 433082294. MS and TR acknowledge support by the Deutsche Forschungsgemeinschaft (DFG, German Research Foundation) - 314838170, GRK 2297 MathCoRe. TW
acknowledges support by the Deutsche Forschungsgemeinschaft 
(DFG) under Germany’s Excellence Strategy within the Cluster of Excellence PhoenixD (EXC 2122, Project ID 390833453).
In addition, we thank Hendrik Fischer for fruitful discussions and comments,
and Fleurianne Bertrand (TU Chemnitz) \auth{and Mary Wheeler (UT Austin)} for some discussions on Mandel's problem.

\bibliographystyle{abbrv}
\bibliography{lit.bib}

\end{document}